\documentstyle[12pt]{article}
\def\={\!=\!}

\def\={\!=\!}

\def\|{{\Vert}}

\def\noi{\noindent}

\def\skip{\vspace{2mm}}
\newcount\u
\def\Remark{\skip\noi{{\bf{Remark \number\u.}}} \advance\u by 1}
 \title{Locally finite knowledge structures}
\author {Robert Samuel Simon}

\begin{document}
\maketitle

\thispagestyle{empty}

\vfill

\noi London School of Economics\newline 
Department of Mathematics\newline 
Houghton Street\newline 
 London WC2A 2AE\newline 
R.S.Simon@lse.ac.uk\newline 
+44-207-955-6753 \newline fax +44-207-955-6877
       
\vfill     
   
\setcounter{page}{-1}
\noi           This research  was supported by 
     the Center for Rationality and Interactive Decision Theory (Jerusalem),  
    the Department of Mathematics of the Hebrew University (Jerusalem),    
 the Edmund Landau Center for Research in Mathematical
Analysis (Jerusalem)  
       sponsored by the Minerva Foundation (Germany),   the German Science 
 Foundation (Deutsche Forschungsgemeinschaft), and a London Mathematical 
 Society Grant to support a visit of Benjy Weiss to London.

\newpage
 \noi 
{\bf Abstract:}
\vskip2cm \noi  In a game of incomplete information, an infinite state space 
 can create problems. When  the  space is uncountably large, the strategy 
 spaces of the players 
 may be unwieldly, resulting in a lack of measurable equilibria.
 When the  knowledge of  a player allows for an infinite number 
 of possibilities, without  conditions 
 on the behavior of the other players, 
  that  player may be unable  to evaluate and compare 
 the payoff consequences of  
 her  actions. 
  We argue that local finiteness is an important  and 
 desirable property, namely that at every 
 point in the state space  
 every player knows that only a finite number of points are possible. 
Local finiteness  implies a kind of  common knowledge of a countable number 
 of points.  
 Unfortunately   its relationship to  other forms of common knowledge is  
   complex. In the context of the  multi-agent propositional calculus, 
 if the set of formulas held in common knowledge is  
   generated 
 by a   finite  set of formulas   but a finite structure is not 
 determined 
  then there are uncountably many  locally finite structures 
 sharing this same set of formulas in common knowledge and likewise 
 uncountably many with uncountable size. This differs 
 radically from the  infinite generation
 of  formulas in common knowledge, and 
 we show some examples of this. 
 One corollary is that  
 if there are infinitely many distinct  points  but 
   a uniform bound on the number of points any  player 
 knows is possible then   the set of formulas in common knowledge cannot 
 be finitely  generated.   \vskip.2cm 
 
\noi    

 \thispagestyle{empty}

    \newpage

\section {Introduction}

\subsection {The problem with infinity} 
 We want to understand  infinite games of incomplete information 
 whose lack of finiteness  lies not in the infinite 
repetition of stages  nor in  the number of actions of the players 
 but rather in the infinite  structure of the
information.   We assume that nature chooses a point in a state space and 
 then the players are informed partially of nature's choice, either 
 through  partitions of the state space or sigma algebras associated to 
the various players. 
    If  the state 
 space is finite, we
remain in the context of a finite game tree,  for which Nash's Theorem
shows that there are mixed equilibria.  The size of the game tree may grow 
  exponentially in the size of the state space, nevertheless it remains 
   finite and 
 subject to the fixed point theorems underlying the Nash Theorem. 
 With an appropriate topology, 
 one can approximate infinite state spaces with   
finite state spaces.
But  how do the equilibria of the approximate finite state space games   
  relate to the equilibria of the original game, if indeed there are any?
\vskip.2cm 
\noi  The space of all functions from the continuum 
to a set of only  two elements is extremely wild, indeed it is equivalent 
 to the set of all subsets of the continuum, a set of higher cardinality 
 than the continuum.  Since the
integral of an arbitrary function from $[0,1]$ to $[0,1]$ is not well
defined, even with one player and an arbitrary strategy 
  we could have trouble defining a
payoff.  
 The  measurability of the strategies is 
an essential issue.  But once we require that  the strategies of the players 
 are   measurable in some sense, the existence of an equilibrium in measurable 
 strategies is  called 
 into question. For fixed point theorems to work
 for non-zero-sum games, usually we need that the
 strategy spaces are compact and convex and that there is  continuity 
  from  the strategies to  the resulting payoffs. By defining 
 the topology weakly we can get compactness of the strategy spaces, however 
 continuity to the payoffs 
 may  fail. By defining the topology strongly we can accomplish the 
  continuity to the payoffs but the compactness of the strategy spaces 
 may fail.    \vskip.2cm 
\noi 
Roughly this is the background to the non-zero-sum  example presented 
 by this author   (Simon 2003), for which 
 equilibria do exist but none of  which are measurable. This example
 demonstrates that approximating by finite state spaces teaches us 
 very little about the original game.  In this example 
 a strong locally finite property holds and local finiteness guarantees that some 
 non-measurable equilibria do exist (Proposition 1, Simon 2003). This focusses 
 attention  onto the local finiteness property as an important condition 
 sufficient  for the existence of equilibria.  \vskip.2cm 
\noi  It should be
 noted  however that   with zero-sum games of incomplete 
 information the situation is quite different.  A weak form of 
 continuity, separable 
 continuity, combined with strategic compactness  suffices for 
 the existence of a value and optimal measurable 
 strategies,  
 proven by Mertens, Sorin and Zamir (1994) with the help 
 of Sion's Theorem (Sion 1958). 
\vskip.2cm 
\noi Next we describe our main results, followed by how they are 
 related to games of incomplete information. 
\subsection{Logic and semantic models} 

 Strictly speaking, our main results pertain to modal logic.
  This perspective  introduces indirectly  a more  
 topological 
 approach. Its basic structure is that of knowledge 
 without any probability attached, however it can be  
   related closely  to probability theory and games of incomplete 
 information, as we will see later. 
\vskip.2cm 

\noi Let $X$ be a  set of primitive propositions, and let $J$ be a 
  set of agents. Although it is legitimate to consider the case of 
either  $X$ or  $J$ infinite, for this paper we will assume
  throughout that both $X$ and $J$ are finite.  
     Construct the set ${\cal L} (X,J)$ of  formulas using  the  
     sets $X$   
     and   $J$ in the
      following  way:\newline
      1) If $x\in X$ then $x\in {\cal L}(X,J)$,\newline
      2) If $g\in {\cal L}(X,J)$ then $(\neg g)\in {\cal L}(X,J)$,\newline
      3) If $g,h\in {\cal L}(X,J)$ then $(g\wedge h) \in 
      {\cal L}(X,J)$,\newline
      4) If $g\in {\cal L}(X,J)$ then $k_j g\in {\cal L}(X,J)$
for every
       $ j\in J$, \newline 
       5) Only formulas constructed through application of the 
        above four  rules are members of ${\cal L}(X,J)$.\newline 
  We write simply ${\cal L}$ if there is no ambiguity. $\neg f$ stands for
 the negation of $f$, $f\wedge g$ stands for  both  $f$ and $g$. 
 $f \vee g$ stands for either $f$ or $g$  (inclusive) and
 $f\rightarrow g$ stands for $\neg f \wedge g$.   
\vskip.2cm 
\noi The connection to games of incomplete information 
 is that an $x\in X$ could represent 
 a fact about the game, for example the validity of a payoff  matrix, the hand 
 that a player holds in a game of cards, or the probability with which a player 
 believes something.

  \vskip.2cm
\noi The relationship between multi-agent logic and games of incomplete 
 information is mediated by semantic models called {\em Kripke Structures}. 
 There are several ways to define Kripke structures, but for 
 our purposes we present the partition model (otherwise 
 known as that corresponding to the S5 logic, which will be 
 explained later). A Kripke structure
 ${\cal K} = ( S, X, J , ({\cal P}^j \ | \ j\in J), \psi)$ 
is defined by two sets $S$ and $X$, 
 a set $N$ of agents, 
 a collection $({\cal P}^n\ | \ n\in N)$ of partitions of $S$ and 
 a function $\psi : S \rightarrow \{ 0,1\} ^X$. The agent $n$ can distinguish 
 between two points in $S$ if and only if they belong to different members 
 of the partition 
 ${\cal P}^n$. The statement 
 $\psi (s)^x=0$ means that $x$ is  not true at $s$ and 
 $\psi(s)^x=1$ means that $x$ is true at $s$.

\vskip.2cm 
\noi We extend the function $\psi$ by defining 
  a map  $ \alpha^{{\cal K}}$ from ${\cal L}(X,J)$ to 
$ 2^S$, the subsets of $S$, 
 inductively on the structure of the formulas, the set $ \alpha^{{\cal K}} (f)$ 
 is where the formula $f$ is valid: 
      \newline {\bf Case 1 $f=x\in X$:}  $ \alpha^{{\cal K}} (x):=  \{
  s\in S\ | \ \psi^x  (s) =1\}$.
     \newline {\bf Case 2 $f=\neg g$:} $ \alpha^{{\cal K}} (f) := 
     S\backslash \alpha^{{\cal K}}(g)$,
     \newline {\bf Case 3 $f=g\wedge h$:} $ \alpha^{{\cal K}} (f) := 
      \alpha ^{{\cal K}} (g) \cap  \alpha^{{\cal K}}
     (h)$,
     \newline {\bf Case 4 $f=k_j(g)$:} $ \alpha^{{\cal K}}(f):= 
     \{ s\ |\ s\in P\in {\cal P}^j 
\Rightarrow P\subseteq \alpha^{{\cal K}} (g)\}$.  
\vskip.2cm 

\noi 
A
 cell of a Kripke Structure  is a member of the
  meet partition 
 $\wedge _{i=1} ^n {\cal P} ^i$, or equivalently a minimal set $C$
 such 
that  for all $j$ the properties 
 $P\in {\cal P}^j$ and $P \cap C \not= \emptyset$ imply  that 
 $P\subseteq C$.  A member of ${\cal P}^j$ for some $j$ is 
 called a {\em possibility} set. 
 A cell has finite fanout if at every point in the cell and   
 every agent $j$ the possibility set 
 of ${\cal P}^j$ containing the point has
 only finitely many elements. Finite fanout is the 
 logical equivalent of local finiteness. 
A formula $f$ is held in common knowledge at a point $x$ of a Kripke 
 structure  if 
 $f$ is true at $x$ and $k_{n_1} k_{n_2} \cdots k_{n_l}f$ is true 
 at $x$ for every choice of a finite sequence $n_1, n_2, \dots n_l$ of  agents.
    A set of formulas in common knowledge is finitely generated 
 if the common knowledge of some finite subset implies the common 
 knowledge of the whole set.

\vskip.2cm 
\noi Given a Kripke structure, 
 construct a sequence ${\cal R}_0, {\cal
 R}_1, \dots $ of partitions of $S$ by ${\cal R}_0 = \{ \psi^{-1} (a)
 \ | \ a\in \{ 0,1\} ^X\}$ 
 and $x$ and $y$ belong  to the same member of ${\cal
 R}_i $ if and only if $x$ and $y$ belong to the same member of ${\cal
 R}_{i-1}$ and for every person $j$ the members $P_x$ and $P_y$ 
 of ${\cal P}^j$
 containing $x$ and $y$ respectively  intersect the same members of
 ${\cal R}_{i-1}$.  Let ${\cal R}_{\infty}$ be the limit of the ${\cal
 R}_i$, namely $x$ and $y$ belong to the same member of ${\cal
 R}_{\infty}$ if and only if $x$ and $y$ belong to the same member of 
 ${\cal R}_i$ for every $i$. 
\vskip.2cm 
 
\noi 
For any set $X$ and set  $J$ of persons there  is a  
canonical Krike structure  
  $ \Omega = \Omega(X,J)$   defined by the formulas in ${\cal L} (X,J)$ 
   such
 that   from any Kripke structure ${\cal K}$ 
  (using the same $X$ and $n$)  there is a canonical map   to
 $\Omega$ defined by the map $\alpha^{\cal K}$  with the property 
 that $y$ and $z$ are mapped to the same point of $\Omega$ if and
 only if $y$ and $z$ share the same member of ${\cal R}_{\infty}$. 
Any Kripke structure with finite fanout 
 will be mapped surjectively  to a cell of $\Omega$ with finite fanout, 
 so that if ${\cal R}_{\infty}$ separates 
 all the points of a Kripke structure  with finite 
 fanout then 
 the structure is isomorphic to a cell of $\Omega$. 
\vskip.2cm

\noi   Our main result is that 
 for every  finitely generated
set of formulas that can be  known in common 
  either this set determines uniquely 
 a  finite cell of $\Omega$   or sharing 
  this same  set of formulas  in common knowledge  there are 
   uncountable many cells of $\Omega$ 
 with  finite fanout and also  uncountably many 
 cells with infinitely many possibility sets 
 of uncountable size.    Furthermore, if there is a uniform bound 
 on the size of the possibility  sets and the cell 
 is infinite   then the formulas in common knowledge 
 cannot be finitely generated.  The situation is very different, 
however, for sets of formulas held in common knowledge that are 
 not finitely generated -- if there are
  uncountably many corresponding 
cells  then either none of these cells or all of them could have  
  finite fanout.  And given  that the  set of formulas held 
 in common knowledge  is not finitely generated and finite 
 fanout doesn't hold    there may be 
 a large difference  between the structure of the 
  cells of $\Omega$  and the cells 
 of some  other Kripke structure holding the same set of formulas 
 in common knowledge.  
\vskip.2cm 

\noi     
\noi  Central to  
 understanding our results   is  
 point-set topology. For every Kripke Structure 
 ${\cal K} =  (S; X;  N; ({\cal P}^n\ | \ n\in N);  \psi)$
  we define a topology
 on the set $S$, the same as in Samet (1990).  
 Let $\{  \alpha^{\cal K}(f)\ |\ f\in {\cal L}\} $ 
 be the base of open sets of $S$. We call this the topology {\em induced 
 by the formulas}.
 The topology of a subset $A$ of $S$ will be the relative  
  topology for which the open sets of $A$ are $\{ A\cap O\ |\ O$ is an open 
   set of $S\} $.       
\vskip.2cm 
\noi To some extent this paper follows results from  two previous papers 
 of this author. In Simon (1999) we showed that  for every set 
 of formulas that can be held in common knowledge either there 
 were uncountably many cells of $\Omega$ corresponding to those 
 formulas or there was only one. The proof used Baire Category 
 in the following way: it was shown that a cell shares its formulas in
 common knowledge 
 with no other cell if and only if  there is some  
  point in the cell and  a positive integer 
  $n$ such that  the set of points that 
 are reachable in $n$ steps from this point 
 (using alternating possibility sets)   
 is an  open set relative to the closure of the cell (meaning that 
 it is not meagre in its closure). 
 A cell that shares its formulas in common knowledge 
 with no other cell is defined to be {\em centered}.   In Simon (2001) 
 we showed that if  the set of formulas held in common knowledge 
  is  finitely generated then the maximality of this set 
 of formulas is equivalent to the finiteness of the corresponding cell 
 and  is equivalent to its being centered. \vskip.2cm

\noi 
Given finite generation of the formulas in common knowledge, 
constructing uncountably many cells with finite fanout 
 is much  more difficult than constructing uncountably 
 many   cells of uncountable size. 
 Even constructing just {\bf one} cell of finite fanout is difficult.  
   To construct a cell of finite fanout dense in $\Omega$ (meaning 
  only the tautologies in common knowledge), 
 we spliced together an infinite sequence of finite Kripke structures (Simon 
 1999). For each non-tautological 
 formula that could be true was associated one of these models and a point
 in that model where this formula was true.  
 The added connections to the other Kripke Structure  
  were distant enough from this point so  that after the  
  splicing each such formula  was still true at its corresponding point.  
  To copy the same approach to construct just one cell 
 of finite fanout with a fixed formula $f$  in 
  common knowledge  one must keep  this formula $f$   
 true at {\bf all} (rather than at just some) 
 points after the  new connections are made.
\vskip.2cm 
\noi There is an additional interest from logic to our main results.
 Finite fanout 
   characterizes the uniqueness of the extension of $\Omega (X,N)$ 
 to the  canonical Kripke structures  using  ordinals beyond 
 the first infinite ordinal $\omega$ (Fagin 1994). Our main results show 
 that  this property of unique extension has little  to do 
 with the set of formulas held in common knowledge.

\subsection{Local and global models} 
 The example of Simon (2003) reveals  that there are at least two levels on 
which  a  game of incomplete information is played, calling for 
 two models of the game, a 
 global model and a local model. It is with the local model 
 that local finiteness (equivalently finite fanout) is important.  
 \vskip.2cm 
\noi The global model is that of measurability. 
  There is a sigma 
algebra 
 ${\cal F}$ on the state  space $\Omega$ 
  for  which  a probability measure  $\mu$ is defined. 
  For  every player $i$ there is a corresponding
  set 
 $A^i$ of actions (either finite or infinite) 
 and  for every choice of $a\in A:=  \prod_{i\in
   N} A^i$ and every  player $i\in N$  there is 
an ${\cal F}$ measurable function 
 $f^i (a) :\Omega \rightarrow {\bf R} $ that represents the payoff to Player
 $i$   if   the players choose 
 the actions $a$ at the point $x\in \Omega$. 
 Corresponding to each player $i$  is a sigma algebra ${\cal
   F}_i$ that is smaller than or equal to 
 ${\cal F}$ (as a collection of subsets of $\Omega$) and 
  represents the private information
 of  Player $i$. For each player $i$ there is   
  a  probability distribution $\mu_i$ defined on 
 ${\cal F}$ (which may or may not be the same as $\mu$).  
A strategy for Player $i$ is  a function $\sigma^i$  from $\Omega$ 
 to the probability distributions 
 of her actions that is measurable with respect to ${\cal F}_i$. 
  Let ${\cal S} _i$ stand for the set of all
 strategies for Player $i$ and let ${\cal S}= \prod _{i\in N} {\cal
   S}_i$
 stand for the set of all
 strategy profiles, a choice of a strategy for each player. 
 \vskip.2cm

\noi  Given that $A$ is finite, for  
  any strategy  profile $\sigma= ( \sigma^i \ |\ i\in N) \in {\cal
  S}$ and every  player $i\in N$   a corresponding payoff ${\cal
  V}^i(\sigma)$ is defined as the  integration according 
 to $\mu_i$  of $f^i$  over $\Omega$ and  
  the strategy profile $\sigma$.
 The strategy profile $\sigma$ is a {\em Harsanyi} equilibrium if 
 for every alternative strategy profile $\hat \sigma$ and $i\in N$ 
$${\cal V}^i(\sigma\ | \ \hat \sigma^i) \leq {\cal V}^i (\sigma)$$
where $(\sigma\ | \ \hat \sigma^i)$ is the strategy profile that is 
 determined by $\sigma$ for all players other than $i$ but determined
 by 
 $\hat \sigma ^i$ for the player $i$. \vskip.2cm

\noi The local model starts from different assumptions, though both models 
 may  have  a  synthesis.  
 For every player $i\in N$ there is a set $T_i$ representing the {\em
   types} of the player $i$. Let $T:= \prod_{i\in N} T_i$ and for
 every 
 $t^i\in T_i$ define the cross section $C_{t^i}$ to be $\{ t^i\} \times \prod_{j\not=i}
 T_j$.  
 For every $t^i\in T_i$ there is a sigma algebra ${\cal F}_{t^i} $ of
 subsets of  $C_{t^i}$  and a  
 probability distribution $\mu_{t^i} $ on $C_{t^i}$ 
  that is defined on   ${\cal F}_{t^i} $.  
 As before $A^i$ stands for the actions of player $i$ (perhaps dependent 
 on the value of $t^i$), 
   and for every choice of $a\in A:=  \prod_{i\in
   N} A^i$ and every  player $i\in N$   there is a payoff  function 
 $f^i (a) :T \rightarrow {\bf R} $ such that its restriction to any
 $C_{t^i}$ is  
  ${\cal F} _{t^i}$ measurable. $f^i(a)$ represents what Player $i$
  receives  at any $x\in T$ if the players choose the actions $a$.  
 A strategy for a player $i\in N$ is any function $\sigma^i $ from $ T_i$ to 
 the probability distributions on its actions, but there may 
 be no measurability requirements.
 A strategy profile $\sigma$ is a 
 choice of a strategy for each player $i$.  \vskip.2cm 

\noi 
In general, a strategy profile may  not define  expected
payoffs  for some fixed $t^i\in T_i$, as there may be problems of 
 measurability. 
However if the actions are finite 
in number and  for any fixed $t^i \in T_i$ the strategies   of 
 the other players when restricted to 
 $C_{t^i}$  are measurable with respect to ${\cal F} _{t^i}$ 
 then one can define for each of her  actions 
 an expected payoff for Player $i$ at the point
 $t^i\in T_i$. 
 A strategy profile is a {\em Bayesian} 
equilibrium for a point $t=
 (t^i \ | \ i \in N) \in
 T$ 
 if  for all $i\in N$ the expected payoffs 
 corresponding 
 to  the strategies of player $i$ can be evaluated in the  corresponding 
 cross section   $C_{t^i}$   and 
    Player $i$ cannot obtain a higher
 payoff according to that  evaluation  by choosing a
 different 
 strategy.  Notice that if the state space has the locally finite property 
  and there are finitely many actions  then such local 
 evaluations are not problematic even if the strategies 
 used by the other players are wild. As stated above, 
 through  approximation by finite 
 games local finiteness does  imply 
 the existence of a Bayesian equilibrium 
 (Proposition 1 of Simon 2003). \vskip.2cm 

\noi In some ways the global model is the stronger one.
 Under natural conditions (such as a Polish state  space) 
  a global model induces a local 
 model through  regular conditional probability distributions (conditional 
 probabilities  for all the  Borel sets done in a consistent manner, see 
 Breiman 1992). 
Likewise under natural conditions 
   an equilibrium of the global model will induce 
 an equilibrium of the local model almost everywhere (see for example
 Proposition 2 of 
 Simon 2003, whose proof idea was explained to me by J.-F. Mertens). 
However a lack of a global Harsanyi equilibrium in 
 the presence of local Bayesian equilibria   shows that 
 the local model has its advantages, especially when combined 
 with local finiteness. \vskip.2cm 
\noi {\bf Question 1:} Let the player set $N$ be finite, let there be 
 finitely many stages and on each stage  each player has 
 finitely many actions. 
 Assume that the state space $T= \prod_{i\in N} T_i$ is compact 
 and Polish with a Borel probability distribution
 $\mu$. Assume that there is a uniform bound on 
 the payoffs and  for every fixed 
 collection of actions the payoffs to each
 player is a Borel measurable function on $T$. Is there such a game 
   that lacks a Harsanyi equilibrium and for every choice of 
  regular conditional probability distributions for the players there are no  
    strategies that  define a Bayesian equilibrium 
 almost everywhere? 
\vskip.2cm 
\subsection{Ergodic Games} 

It is worthwhile to look briefly at the structure behind the Simon (2003) 
 example. This example belongs to the category of {\em ergodic games}, 
 which will be defined below. 
\vskip.2cm 
\noi 
 The Bernoulli  space can be represented 
 as   the set $\{ a,b\} ^{\bf Z}$, where 
 ${\bf Z}$ stands for the integers (including the negative integers) and 
 $a$ and $b$ are distinct symbols.  
  A point of $\{ a,b\} ^{\bf Z}$ is a doubly infinite sequence  ${\bf x}=
 ( \dots , x^{-1} , x^0 , x^1 , \dots )$ where for every $i\in {\bf
   Z}$ 
 $x^i$ is either $a$ or $b$. The probability distribution on the state 
 space is the canonical one that gives every choice of $a$ or $b$ at 
 any position equal weight and   independently of the choices 
  for  $a$ or $b$ at the other positions.
 The function $T: \{ a,b\} ^{\bf Z} \rightarrow \{ 0,1\} ^{\bf Z}
$, called the shift operator, is defined by 
 $T(x) ^i = x^{i-1}$.  

\vskip.2cm 
\noi In the Simon (2003) example there are two 
 sets of payoff matrices corresponding to $a$ and $b$ and 
  the $0$-coordinate of a point in $\Omega$ 
 determines the payoff matrices, 
 so that if $y\in \Omega$ and $y^0=a$ then payoffs at the point 
 $y$ are determined 
 by the matrices corresponding to   $a$.
There are three players, Players One, Two, and Three.  
Let  $\sigma: \Omega 
\rightarrow \Omega$ be the measure preserving involution defined by   
 $(\sigma (y))^i: = y^{-i}$, where  
 $y^i$ is the $i$th coordinate of $y\in \Omega$. $\sigma$ is the reflection 
 of the  sequence about the position zero.
 Let   $\tau :\Omega \rightarrow \Omega$ be the measure preserving 
 involution   defined by 
$(\tau (y))^i:= y^{1-i}$.  It follows that 
 $T:=\tau \circ  \sigma  $ is
 the usual  Bernoulli shift operator $(T (y))^i = 
 y^{i-1}$. Define the beliefs $t^1$,  $t^2$ and $t^3$ of the players 
   according to   $\sigma$ for
 Player One and $\tau$ for the other players;  this means that 
 at any point $y\in \Omega$   Player One considers only 
  $y$ and $\sigma (y)$ to be possible, and with equal probability;
 (if $\sigma (y) =y$ then 
 Player One believes in $y$ with full probability).   At any $y$, both 
Player Two and Player Three believe that only $y$ and $\tau (y)$ 
 are possible, and with 
equal probability.  
 Player Two and Three have the same beliefs.  Player One 
always knows the payoff 
  but not always what the other players 
might know. \vskip.2cm 
\noi The secret to this example is that the shift operator $T$ is an ergodic 
 operator that  acts almost 
 everywhere 
 upon any equilibrium of the game in a way so that this equilibrium 
 cannot be measurable.  
\vskip.2cm 

 \noi  The game example  in Simon (2003) is also an {\em ergodic game}. 
 An ergodic game  is a  game 
 satisfying the following properties: \vskip.2cm 
\noi
(1) there is  one stage of play with moves chosen simultaneously, \vskip.2cm 

\noi (2) there are  finitely many 
 players,\vskip.2cm 
\noi  
(3) each player has  finitely 
 many moves, \vskip.2cm 
\noi  
(4) there is a compact  Polish   space $\Omega$ with an atomless 
  Borel probability distribution 
 $\mu$ representing a choice by nature, 
 \vskip.2cm 
\noi (5) for every combination of moves, one for each player, the
 payoff to each player is a continuous function from $\Omega$ to ${\bf
 R}$, \vskip.2cm      
\noi (6) at every point $x\in \Omega$ every player $j$ has a
 discrete probability
 distribution on $\Omega$, called his {\em belief}, 
 with a finite support set $S^j(x)$   
   containing $x$  such that 
 at all the other points in  this finite support set $S^j(x)$ the player 
 $j$ has the same discrete distribution, \vskip.2cm 
\noi   
(7) these discrete 
 beliefs of the players
   change continuously (with respect to 
 the $\mbox {weak}^*$ topology), and for any player $j$ 
 it forms a  regular conditional probability
    of  $\mu$ with respect to the sigma algebra 
${\cal F}^j:= 
\{ B\ | \ B$ is Borel and $x\in B\Leftrightarrow S^j(x)\subseteq B\}$, 
\vskip.2cm 
\noi (8) If $B$ is a Borel set  with $B\in {\cal F} ^j$ for all
 players $j$ then $\mu (B)$ is equal to either zero or one. \vskip.2cm
\noi By Property 6 all ergodic games have the local finiteness property.
 \vskip.2cm

\noi  
 A measure preserving transformation $T$ on a probability space 
 with a Borel probability distribution $\mu$ is
called {\em ergodic} if the only Borel sets $B$ with $T^{-1} (B) = B$
satisfy 
$\mu (B) =0$ or $\mu (B) =1$.  
  Why is the term ergodic used to describe these games?  The ergodic 
 aspect of these games is contained in Property 8. It 
 guarantees that the game cannot  not 
 decompose into two different subgames of positive  probability. 
\vskip.2cm 

\noi 
 There is an easy  way to create an ergodic game with $n$ players.
 Let $\Omega$ be a compact  Polish probability space with a Borel 
  probability measure $\mu$. 
 A function $f: S \rightarrow S$  is called an {\em  involution}
if $f$ is not the identity function but $f^2$ is the identity function.  
  Let $\sigma_1, \sigma_2, \dots ,
 \sigma_n$ be a set of $n$ continuous 
 measure preserving involutions such that
 some combination of the $\sigma_i$'s is an ergodic transformation.  
For each player $i$ define ${\cal F}_i$ to be $\{ B \ | \ B \mbox {
  Borel }, x\in B \Rightarrow \sigma_i (x) \in B\}$ and $T_i$ to be
the equivalence classes defined by $x \sim_i  y$ if and only if 
 $x=y$ or $y= \sigma_i (x)$. 
For every player $i$ define
 the belief  function $t_i:\Omega \rightarrow \Delta (\Omega)$ by  
 \newline $t(x) (B)=$ \newline 
$0$ if both $x$ and $\sigma_i (x)$ are not in $B$, \newline 
$1$ if both $x$ and $\sigma_i (x)$ are in $B$, and \newline 
$1/2$ otherwise. \newline 
 $t_i$ is  continuous with respect to  the $\mbox {weak}^*$
 topology  and 
 the function $t: \Omega \rightarrow \Delta (\Omega)$ 
 is a regular conditional probability  induced by  
$\mu$ and ${\cal F}_i $ (Lemma 0 of Simon 2003). 
The rest of the properties of an ergodic game can be constructed  easily.
 Indeed the same can be done for any collection of finite groups $(G_i\ | \ 
 i\in N)$ that act measure preserving on a Polish probability space. In this 
 case at any point $x$ the player $i$ would consider the points 
 of the orbit 
 $G_i x$ to be possible with each point $y= gx$ weighted according 
 to the number of group elements $g\in G_i$ satisfying $gx=y$.

\subsection{Absolute Continuity}
 There is a  condition on a game of incomplete information  
  which, along with natural assumptions on the payoffs, 
  guarantees the existence of measurable equilibria (as demonstrated 
 by Milgrom and Weber 1985), {\em absolute 
 continuity with respect to the marginals}. With a weak topology on 
 the strategies, it implies  
  continuity from a compact  strategy space to the payoffs. 
  It is worthwhile 
 to see how this condition may fail, especially when local finiteness  holds.

   \vskip.2cm 
\noi  If $S = A \times B$  and $\mu $ is a probability
 distribution on $S$ then the marginal probability 
distribution on $A$  gives 
 $\mu (C \times B) $ to any measurable subset $C$ of $A$ (and likewise
 define the marginal probability distribution on $B$).   A measure 
 $\mu$ is {\em absolutely continuous} with respect to another measure 
 $m$ if and only if $m (A) =0$ implies that $\mu (A)=0$.  One way for
 a measure $\mu$ to be absolutely continuous with respect to another
 measure $m$ is for $\mu (B) = \int f(x) 1_{B} d(m) $ for some
 non-negative and integrable function 
 $f$. The Radon-Nikodym Theorem states that in many common  situations the
 converse holds, namely that absolute continuity implies the existence
 of such a function. \vskip.2cm  
\noi  A game of incomplete information  with a probability 
 measure $\mu$ on a state space $\Omega = \prod_{j\in
 N} T_j$  (a synthesis of the global and local models) 
 has the  {\em absolutely
 continuous}  property  if $\mu$ is absolutely continuous with 
 respect to the product of the marginals induced by $\mu$  
 on  each of the $T_i$. 
\vskip.2cm 
\noi Notice that no ergodic game can have the absolute continuity property, 
 indeed neither can any 
  game with the local finiteness property where the Radon-Nikodym 
 and Fubini theorems hold and every singleton set has zero measure. 
  It suffices to prove this  for two players, as the proof
 for more players introduces no new ideas. 
For each  player $i\in \{1,2\}$  let 
  $\mu_i$ be the marginal distribution of $\mu$ on $T_i$. As
 any single point of $\Omega$ is given zero probability by $\mu$ and
 every cross section of $\Omega$ 
 defined by  a point of $T_i$  contains only finitely many points it
 must hold that every single point of $T_i$ is given zero probability
 by $\mu_i$.  By Fubini's Theorem we can re-write  $\mu (\Omega) = 
 \int h(x_1, x_2)\  d (\mu_1 \times \mu_2)  $ as 
 $\int ( \int  h(x_1, x_2)  d \mu_1 (x_1))\  d\mu_2 (x_2)$, 
where $h$ is the function implied by the  
 Radon-Nikodym Theorem. 
 But for
 every fixed $x_2\in T_2$ there are only finitely many points $x_1$ in 
 $T_1$ with $(x_1, x_2) \in \Omega$.
 Because $\mu_1$  gives zero probability to every single
 point it must follow that $\int  h(x_1, x_2)  d \mu_1 (x_1)$ is 
 zero for every fixed choice of $x_2$, hence $\mu (\Omega)=0$,
 a contradiction.       
\vskip.2cm 
\noi 
Any game where two distinct players always  share identical  
 information and the marginals are atomless 
 will also  fail the absolute continuity 
 property, because any set that projects canonically 
 onto  the diagonal of the product of these two players' type spaces 
 will be a set of measure zero in the product topology. 
\vskip.2cm 
\noi Consider also the following information structure on which a game 
 can be defined for which  local finiteness holds but 
 the absolute continuity   
 property fails. 
Let  $\alpha$ be an irrational real number with  
 $0< \alpha <1$ and let $C$ be the  subset of the square $ [0,1]\times [0,1]$ 
 defined by $C= \{ (x,x) \} \cup \{ (x,y) \ | \ y=x+ \alpha$ 
 or $y=x+\alpha -1\}$.  
  Let    $\mu$ be the  probability  measure on $C$ defined by   
 \newline 
 \noi (1) $\mu   (\{ (x,x) \ | \ x\in [a,b]\}) = \frac 12 (b-a)$, 
  \newline  
\noi (2) $\mu (
\{ (x,y) \ | \ y=x+ \alpha$ 
 or $y=x+\alpha -1, x\in [a,b]\}) = \frac 12 (b-a)$. \vskip.2cm 
\noi  To make this example into a two-player 
game one can  define payoff matrices 
 that change continuously according to the location on $C$. Due to 
 the special locations where the cross sections will be three 
 points instead of the usual two, strictly 
 speaking such a  game will not be ergodic. However by identifying 
 the values of $0$ and $1$ for both players an ergodic game can be 
 constructed.   It would be interesting to discover whether a non-zero-sum 
 game 
 can be so constructed where no measurable equilibria exists. 
\subsection{Finite additivity}

\noi One  could think  
 that the countably additive axiom    of the conventional 
 definition of 
 a probability distribution is to blame for the discrepancy between 
  local and global equilibria, in particular   the possibility of 
  Bayesian 
local equilibria where there is no global Harsanyi 
  equilibria. (It should be noticed, however, 
 that  the proof that a global Harsanyi equilibrium induces 
 a local Bayesian equilibrium almost everywhere does use countable additivity.)
 Indeed on the Bernoulli shift space 
 if one required only finite additively there will be 
  many finitely additive shift-invariant measures  
  that are defined on all subsets. This is due 
 to the fact that the shift transformation defines an amenable group action. 
 However with  $G_i$ subgroups defining the beliefs 
 of the players (as described  above) 
 generating  a non-amenable group action 
 there will be  
  state spaces  for which there are no 
 finitely additive measures defined on all the subsets 
 that are also invariant with regard to this group.
 If this  groups acts 
 on any measurable equilibria in an appropriate way it may be 
 possible to demonstrate the non-existence of finitely additive  
 global Harsanyi equilibria, though there will  be many local Bayesian 
 equilibria.    \vskip.2cm 
\noi 
 In this context it may be  relevant to review the related   
   Banach-Tarsky Paradox  from 1924, which speaks directly to the 
 non-existence of finitely additive measures defined on all subsets  that are 
 rotation invariant. 
 The paradox states that there is a way to partition the sphere $S^2$  into
finitely many parts $A_1, \dots , A_k , B_1, \dots, B_l$
  such that after rotating these parts  two copies of $S^2$ 
are created, one sphere from the $A_1, \dots , A_k$ and another sphere
from the $B_1, \dots , B_l$. The group of rotations of $S^2$ is 
 a non-amenable group. To review the relations between amenability, 
  the paradox, and other  issues, see Wagon (1985). \vskip.2cm  
\noi With inspiration from the Banach-Tarsky paradox and its  
 relation to ergodic games through amenability,
 we pose the following 
  related open questions. 
\vskip.2cm 
\noi {\bf  Question 2:} Does there exist an ergodic game that
has no equilibria measurable with respect to any finitely additive
probability measure of the state space?
\vskip.2cm \noi  The analogy between Question 2  
  and the Banach-Tarsky paradox is the following.  Assume the information structure of 
 the ergodic game is generated through 
 finite groups representing the beliefs of the players as described above. 
  The subsets $A_1, \dots, A_k,B_1, \dots , B_l$, dependent on 
 the strategies used,  could 
cover the  state  space  and   
 represent the locations where certain  subsets of strategies 
 are used by the players. One would like  to show 
 that in any Bayesian equilibrium there will be group elements acting on 
 these sets in a way similar to the rotations in  the Banach-Tarsky paradox, 
 demonstrating that the equilibrium cannot be measurable with respect 
 to any finitely additive measure. There is a parallel in the main proof 
 of Simon (2003), where it is shown that if some 
 player alternates her  behavior 
 throughout the state space between
 two pure strategies and the measurable  subsets $A,B$ 
 represent the locations where these two pure strategies are used then 
 the ergodicity of the square $T^2$ of the 
  shift operator $T$   and the property $T(A)=B$ and $T(B)=A$ almost  
 everywhere implies that either both $A$ and 
 $B$ are of measure zero or both $A$ and $B$ are of measure 
 one, both contradictions. \vskip.2cm 
 
\noi A question related indirectly to Question 2 is the following: \vskip.2cm
 \noi {\bf  Question 3:} Does there exist an ergodic game
 and  some positive $\epsilon>0$ such that the game does
 not have a global Harsanyi  $\epsilon$-equilibrium in Borel  measurable 
 strategies? 
\vskip.2cm 
\noi The example  in Simon (2003) relies 
 heavily on the non-linear aspects of how the payoff can change 
 when one player fixes her  action and the 
 other two players vary their mixed strategies. This inspires the following 
 question, which can be amended to refer to finite additivity and measurable 
 $\epsilon$-equilibria: \vskip.2cm 
\noi  
 {\bf Question 4:} Is there a two-person non-zero-sum 
 ergodic game  that has no global  Borel measurable  equilibrium?    
\vskip.2cm 

\noi One 
  could  think  that  local finiteness 
 is not  appropriate  for the  context of 
 equilibrium existence, rather   
 the locally countable  property,   namely that every player at every 
 point knows that one of countably many points are possible. 
 Indeed  a countably additive 
  measure  on a countable set  will be determined by the weights 
 given to all the individual points, and indeed 
 the proof of local Bayesian equilibria in Simon (2003) is 
extended to such a context.
 But there are 
 limitation to local countability that are not present in local 
finiteness. For  finitely  additive measures 
 there is no such unique determination, as there can be many distinct 
 finitely additive measures on the same countable set that assign 
 zero to all singleton sets. Furthermore 
 the infinity of any set introduces topological questions that do not arise 
 with finite sets, namely which sequences of points  converge, and if so 
 will this convergence point be a member of that same set. This topological 
  complexity 
 manifests itself further   when we consider 
 the logical aspects.

\subsection{Meagre and null sets}

Although measure theory and Baire Category are distinct approachs, 
 they have many parallels.  In most contexts 
 a  meagre or no-where dense 
 set will have zero measure, so that a first Baire Category set, the 
 countable union of meagre sets, will also have zero measure. 

\vskip.2cm 
\noi More  parallels can be  found  when comparing Kripke structures to  
  state spaces of games of incomplete information. 
If $C$ is a cell of $\Omega (X,N)$ then 
 the  formulas in common knowledge in $C$ define  the  closure  
 of the  cell, parallel  to the definition of the {\em support} 
 of a probability measure on a compact and separable probability  
 space  as the smallest closed set of probability 
 one.  Property   8 of ergodic games, that a set known in common by 
 the players must be either measure zero or measure one, finds its parallel 
 in the Baire Category argument underlying the centered vs uncentered 
 distinction of the cells of $\Omega (X,N)$.  
\vskip.2cm 
\noi It should be noted  that it is easy to construct a probability
 distribution 
 for a Kripke structure,  forming the basis of   
  a game of incomplete information. Consistent determinations of 
probabilities for all the formulas   will induce a Borel probability
 distribution by way of the Kolmogorov extension theorem (where Borel 
 refers to the topology defined by the formulas).  Such consistent 
determinations start with the probabilities for the validity  of each 
 of the  $x\in X$ and then the probabilities that  the players know or 
 don't know  
 any of the $x\in X$ to be valid.  
 Because this  distribution will be unique on the Borel 
 sets, given that there is Hausdorff separation there will 
 be no need to consider anything but the Kripke structure's canonical 
 map into $\Omega$,  
  unless one wishes to extend this  measure beyond the (null set) completion 
 of the  Borel sets  relative to this probability measure. (Because the 
 base of open sets are defined by the formulas and they define both open and 
 closed sets, any topological separation of points will be Hausdorff.)
  Given that the  
  topology induced by the formulas 
 is not Hausdorff (meaning that the map into $\Omega$ is 
 not injective), one could accomplish separation by 
 introducing more elements to the set $X$. In this context it may 
 noticed that  $X$ can be extended to a  countable set without altering 
 the Baire Category  structure to the cells (Simon 1999).   
 \vskip.2cm 

\noi It is plausible, with a probability distribution  constructed for  
 a Kripke structure and the main results
 of this paper, that one could answer Question 1 
 by   
 demonstrating  a game of incomplete information where there are no 
 measurable Harsanyi equilibria and Bayesian equilibria exist for 
 only an unmeasurable subset or a subset of measure strictly 
  less than one. \vskip.2cm 
 
\noi Going the other way, from either a local or global model for a game 
 of incomplete information to a corresponding Kripke structure, it is easy 
 as long as there are well defined supports for  the local beliefs 
 of  the players (defined directly or indirectly by 
 the model).   \vskip.2cm 
\noi For example, consider the Bernoulli shift space $ \{ a,b\}^{\bf Z}$ 
 where $X= \{ a, b\}$ and $a$  is true at $\omega \in  \{ a,b\}^{\bf Z}$ 
if and only if 
 $\omega^0=a$ (and otherwise $b$ is true), and let 
 the knowledge of Players One and Two be that as defined 
 above by the involutions  in  the example of  Simon (2003),
  and for now we drop the third player who has 
 identical information to the second player.
 Because the shift operator 
 is related to the operations $k_1$ and $ k_2$ on the formulas, 
 one can show easily that  
  the Bernoulli shift space as a Kripke structure 
 will map isomorphically   to a 
 compact subset 
  of $\Omega (X, \{ 1,2\})$  With 
probability one, any point will be in a dense cell and due to  
 local finiteness  there must be  uncountably many dense cells. Because 
 of some 
 exceptional points in the shift space, such as $(\dots, a, a,a, \dots)$, 
 there will be some cells that are not dense. However by adding 
 a new third  player associated with the involution $\pi$ defined by 
 $\pi (z) ^i= z^i$ if $i\not= 0$ and otherwise $\pi$ switches the $a$ with 
 the $b$ and vice-versa  on the $0$-position, one gets a compact subset 
 of $\Omega (X, \{ 1,2,3\})$ that is comprised entirely of uncountably many 
 dense cells.  
 \vskip.2cm

\noi 
The rest of the paper is organized as follows.  Section 2 provides background 
 information. Sections 3 and 4 contain the  proofs of the two 
 main claims, Theorem 1 and Theorem 2. The last section discusses the lack 
 of finite generation in more detail. 
\section {Background}

       \subsection {Formulas and logic}

\noi Recall the above definition 
 of the formulas ${\cal L}$. 
 There is a very elementary  logic defined on the formulas in ${\cal
         L}$  called $S5$. 
For a longer  discussion of the $S5$ logic, see   
 Cresswell and Hughes (1968); 
 and for the multi-person  variation, 
  see Halpern and Moses (1992)       
  and also Bacharach, et al, (1997).
      Briefly, the $S5$ logic  
    is defined by two rules of inference, modus ponens and necessitation,  
    and  five types of axioms. 
      Modus ponens means that if $f$ is a theorem and 
       $f\rightarrow g$ is a theorem, then $g$ is also a theorem.   
       Necessitation 
      means that if $f$ is a theorem then $k_j f$ is also a theorem 
       for all $ j\in J$. The axioms are the following, for every   
        $f,g\in {\cal L}(X,J)$ and $ j\in J$:  
   \newline 1) all formulas resulting from theorems of the propositional  
    calculus through substitution, 
   \newline 2) $(k_j f\wedge k_j (f\rightarrow g) ) \rightarrow k_j g$,
   \newline 
     3) $k_j f \rightarrow f$, 
     \newline 4) $k_j f\rightarrow k_j (k_j f)$, 
     \newline 5) $\neg k_j f\rightarrow k_j (\neg k_j f)$. 
  \vskip.2cm 
  
\noi 
    A set of formulas  ${\cal A}\subseteq 
    {\cal L} (X, J)$  is called  {\em complete} 
    if for every
       formula $f\in {\cal L} (X,J)$  either $f\in {\cal A}$ 
       or $\neg f\in {\cal A}$.
      A set of
      formulas is called {\em consistent}
      if no finite subset
 of this set leads to a logical contradiction, meaning 
  a deduction of $f$ and  $\neg f$ for some formula $f$.
We define $$ \Omega (X,J):=\{ S\subseteq {\cal
L}(X,J)\ |\ 
 S\mbox{ is complete and consistent} \}.$$  A {\em tautology} of ${\cal L} (X,J)$  
 is a formula that is true at every point of $\Omega (X,J)$, or 
 in other words a formula $f$ is a tautology if and only 
 if for every $z\in \Omega (X,J)$ the formula $f$ is in $z$. \vskip.2cm 
   
\noi  The 
$\Omega (X,J)$ is itself a  Kripke Structure
 $(\Omega (X,J); X; J; ({\cal Q}^j(X,J)\ | \ j\in J);  
\overline \psi (X,J))$ with  for every  $j\in J$ the    partition  
      ${\cal Q}^j(X,J)$ being  that 
  generated by the inverse images of  
 the function
  $\beta ^j: \Omega(X,J) \rightarrow 2^{{\cal L} (X,J)}$  defined by 
 $$\beta ^j (z):= \{ f\in {\cal L} (X,J)\ | \ 
  k_j f \in z\}, $$ the set of formulas known by person $j$ and 
 $\overline \psi$ defined by 
 $\overline \psi ^x(z)= 1$ if and only if $x\in z$.       
  Due to the fifth set of axioms $\beta ^j (z)\subseteq \beta ^j (z')$ 
   implies that $\beta ^j (z) =\beta ^j (z')$. 
  We will write $\Omega$, ${\cal L}$, $\overline \psi$  and ${\cal Q}^j$ if there is 
   no ambiguity.

\subsection{Common knowledge}  

      \noi The expression  $E(f)=E^1(f)$ 
    is defined to be $\wedge_{j\in {J}} k_j f$, 
       $E^0(f):=f$, and 
       for $i\geq 1$, $ E^i(f) := E(E^{i-1}(f))$. 
   A formula $f\in {\cal L}(X,J)$ is common knowledge in a subset of 
    formulas $A\subseteq {\cal L}(X,J)$ if 
    $E^nf\in A$ for every $n<\infty$. 
      \vskip.2cm

\noi  The set of  formulas that are held in common
 knowledge is constant within any cell and equal to all the formulas 
 true everywhere in that cell (Halpern and Moses 1992). In
 other words: \vskip.2cm 

\noi   {\bf Lemma A:} 
          For any cell  $C$ of $\Omega (X,J)$
         $\{ f \in  {\cal L} (X,J)\ | \ f $ is common knowledge in $z $
          for some $z \in C\} = \{ f\in  {\cal L} (X,J)\ |\ 
          f $ is common knowledge
          in $z$ for all $z \in C\} = 
          \{ f\in  {\cal L} (X,J)\ |\  f\in z $ for all $ z \in C\} $. 
          \vskip.2cm

     \noi  Due to Lemma A, we have a map $F$ from the set 
      ${\cal Q} = \wedge _{j\in J} {\cal Q}^j$
 of cells  to  subsets of formulas defined by             
      $F(C):= \{ f\ |\ f $ is common knowledge in any (equivalently all) 
      members of $C\}$.

\vskip.2cm

     \noi 
 For any subset of formulas $T\subseteq {\cal L}$ 
 define  
 $\underline {Ck}(T):=\{ f\in {\cal L}\ |\ $  
  there exists 
   an $i<\infty$ and a finite set $T'\subseteq T$   
   with 
 $(\wedge_{t\in T'} E^i(t)) \rightarrow f$ a tautology $\}$. 
         We 
         define ${\cal T}(X,J)=\{ \underline {Ck}(T)\ |\ T\subseteq 
         {\cal L}(X,J)\}\backslash $   $\{ {\cal L}(X,J) \}, $ 
        and  we say that $T$ 
     generates $\underline { Ck}(T)$.  
         If there is no  
          ambiguity, we can write simply ${\cal T}$. 
        $\underline {Ck}(T)$ is the set of formulas whose common knowledge 
         is implied  by the common knowledge of the formulas 
          in $T$.
\vskip.2cm

 \noi 
 An $S\in {\cal T}$ is finitely generated if there exists a finite subset 
 $T\subseteq S$ such that $\underline {Ck} (T)=S$. 
For every set of formulas $T\subseteq {\cal L}$ define  
 the set $${\bf Ck}(T):=\{ z \in \Omega\ |\ \mbox { every member of } 
 T\mbox{ is 
  common knowledge in } z\}.$$ 
      For any $T\subseteq {\cal L}$,  
      ${\bf Ck}(T)$ is a closed set, since the     
       ${\bf Ck}(T)$ is the intersection of the sets   
      $\alpha(E^l f)$ for all $l<\infty$  and  all 
       formulas $f$  in $T$. \vskip.2cm 

\noi   
A cell $C$ is defined to be {\em centered} if and only if there is 
 no other cell $C'$ with $F(C') =F(C)$. 
     \vskip.2cm 
\noi  
In Simon (1999) we proved that if a cell $C$ is not centered then 
 there are uncountably many other cells $C'$ such that $F(C')=F(C)$. 
Since $\underline {Ck}(T)=F(C)$ means that the cell $C$ is dense 
 in ${\bf Ck}(T)$; the cell $C$ being not centered is equivalent to the 
 existence of 
uncountably many other cells $C'$ that are also dense 
 in ${\bf Ck}(F(C))={\bf Ck}(T)$.

\vskip.2cm

\noi 
 In this paper we prove that  if $C$ is a cell and  $F(C)$ is  
  finitely generated then either $C$ is finite and there is no 
 other cell $C'$ with $F(C)=F(C')$ {\bf or} 
   there is a continuum of   distinct cells $C'$ with continuum 
 cardinality such that $F(C')=C$  and there is a continuum  of  
 distinct cells $C'$ of finite fanout such that  $F(C')=C$.

\subsection {Kripke Structures} 

    \vskip.2cm 
\noi   If ${\cal K} = (S;X; J; ({\cal P}^j\ | \ j\in J); \psi)$
 is a Kripke structure  we  
      define a map $\phi^{{\cal K}}:S\rightarrow \Omega (X,J)$ 
       by
      $$\phi^{{\cal K}}(s):= \{ f \in {\cal L}(X,J)\ | \  \ 
     s\in  \alpha^{{\cal K}} (f)\},$$ 
 where $\alpha ^{\cal K}$ is the map defined above.  
This is the canonical map, also contained
     in Fagin, Halpern, and Vardi (1991). 
\vskip.2cm

\noi  
 {\bf Theorem:} For every $f\in {\cal L}(X,J)$, $f$ is a theorem of 
  the multi-agent S5 logic if 
  and only if $f$ is a tautology. 
  Furthermore, 
    $\phi^{\Omega} (z) = z $ for every $z\in \Omega$.
          \vskip.2cm
          
     \noi For a proof of the first part of this theorem, see Halpern and 
     Moses (1992) and  Cresswell and Hughes 
     (1968), and for how the second part follows from the first part  see 
     Aumann (1999).  
      We will call this result the ``Completeness Theorem."

 \vskip.2cm

           \noi  For a Kripke Structure  
           $  {\cal K} =(S;X; J; ( {\cal P}^j \ | \ j\in J) ;  \psi )$,   
          if $s\in \alpha^{{\cal K}} (f)$, or equivalently  
           $f\in \phi^{{\cal K}}(s)$, we say that $f$ is true  
           at $s$  with respect to ${\cal K}$. 
           We say that $f$ is valid  with respect to 
           the Kripke Structure ${\cal K}$ 
            if $f$ is true at $s$ with respect to
 ${\cal K}$ for every $s\in S$. 
The Kripke Structure is                 
          {\em connected} if the meet 
          partition $\wedge _{j\in 
          J}{\cal P}^j$ is a singleton (equal 
           to $\{ S\}$).  We define a {\em connected component} of 
            a Kripke Structure to be a member of this  
          meet partition.
          Two points $s,s'\in S$ are {\em adjacent} if they share some 
           member of ${\cal P}^j$ for some $j\in J$. 
          We 
      define the {\em adjacency distance} 
     between any two 
      points $s$ and $s'$ in $S$ as  
      $  \min \ \{ d\ |\ $ there is a sequence $s=s_0 , \dots 
     , s_d=s'$, a function $a : \{ 1,\dots , d\} \rightarrow J$ 
      and sequence of sets $D_1\in {\cal P}^{a(1)},\dots, 
      D_d\in {\cal P}^{a(d)}$ such that 
       for all $ 1\leq i\leq d\quad s_i$ and $s_{i-1}$ 
both belong to $D_i\} $, 
      with zero distance between any point and itself and infinite distance  
      if there is no  
      such sequence from $s$ to $s'$. 
Such a sequence we call an {\em adjacency path}.
\vskip.2cm 

\noi 
Given a 
 Kripke Structure  ${\cal K}=(S; X; J ;
({\cal P}^j\ |\ j \in J); \psi )$ and a subset $A\subseteq S$,  
          we define another Kripke Structure 
 ${\cal V}^{\cal K}(A):= 
 (A;X; J ; 
({\cal P}^j|_A  \ |\ j \in J)); \psi|_A )$ 
 where for all $j\in J$ 
 ${\cal P}^j|_A:= \{ F\cap A\ |\ F\cap A\not= \emptyset $ and $F\in  
  {\cal P}^j\} $  and for all $x\in X$ and $a\in A$ 
   $\psi^x|_A(a)=\psi^x(a)$. 
   If there is no ambiguity concerning the initial model  
    $\cal K$, we can replace ${\cal V}^{\cal K}( A)$ by 
   ${\cal V}(A)$.

\vskip.2cm 
\noi 
Now we can see why a Kripke Structure 
 with the finite fanout property is 
 essentially a cell with finite fanout. 
 It is easy to prove that for every  set $P \in {\cal P}^j$
 corresponding to  an agent $j\in J$ 
 in a Kripke Structure 
 ${\cal K}$ that $\phi^{\cal K}(P)$ is a dense 
 subset of some member of ${\cal Q}^j$ (Lemma 5, Simon 1999).  
 Fagin (1994) proved that a cell has a unique extension to 
 all canonical Kripke Structure corresponding to the transfinite ordinal 
 numbers beyond the first infinite ordinal if and only if it is of finite 
 fanout, and that  representation in all these canonical Kripke Structure 
 characterizes  the interactive knowledge of the agents.   Combining 
 these results, restricting ourselves to $\Omega$ is sufficient for 
 understanding  Kripke Structures with finite fanout.  

    \subsection {Canonical Finite Models}

Behind  our main results is a hierarchical 
 construction of $\Omega$.  Every formula has a ``depth", an inductively 
defined natural number representing the extent to which   the knowledge 
 operators $k_j$ of the agents $j\in J$ 
have been  used to define the formula.  (Formulas of 
 depth zero are those of the conventional propositional calculus, 
  constructed without the knowledge operators of 
 any agents.)  For every natural number $i$ there is a finite 
 Kripke Structure $\Omega_i$ that represents the knowledge 
 of the agents up to the depth of $i$. Furthermore   
  $\Omega$ is the inverse limit of the $\Omega_i$, meaning that a point 
 in $\Omega$ is defined by a sequence of extensions 
 from $\Omega_i$ to $\Omega_{i+1}$ for all the $i$. If the 
set of formulas held in common knowledge is finitely generated, 
 then   these formulas have a maximal depth $d$.  By exploiting the  
 choices in how one could  extend a point in $\Omega_i$ to 
 $\Omega_{i+1}$ for some of the $i$ that are greater than $d$, one can 
 construct cells in the limit of the process that do or 
 do not have finite fanout.  If the corresponding set of formulas 
 held in common knowledge is not finitely generated, the lack 
 of a maximal depth for a generating subset renders the hierarchical    
     construction meaningless for our purposes. 
\vskip.2cm 
             \noi   
  We define the {\em depth} of a formula inductively on the structure of
   the formulas.   If $x\in X$, then depth $(x):=0$. 
  If $f=\neg g$ then depth $(f):=\mbox {depth }(g)$; if $f=g\wedge h$
   then depth ($f$):= max (depth ($g$) , depth ($h$)); and
   if $f=k_j (g)$ then
    depth ($f$):= depth ($g$) +1. 
  \vskip.2cm 
          
      \noi 
       For every $0\leq i<\infty$ we
 define ${\cal L}_i:= \{ f\in {\cal L}\ |\ $depth $(f)\leq i\}$ and 
 define  $\Omega_i$    to be the set of maximally consistent subsets of 
  ${\cal L}_i$.
            If there may be  ambiguity, 
              we will write $\Omega_i (X,J)$.
             We must  perceive an $\Omega_i$ in 
         two ways,  
          as a  Kripke Structure in its own right and  
          as a canonical projective image of $\Omega$ inducing 
           a partition of $\Omega$ through inverse images.
          We define   
            $\pi_i:\Omega \rightarrow \Omega_i$ to be the  
           canonical projection $ \pi_i (z) := 
  z \cap {\cal L}_i. $   
            Due to an application of Lindenbaum's Lemma,  
            the maps $\pi_i$ are surjective. For all $i\geq k$ define the map 
 $\pi^i_k$ to be $\pi_k \circ \pi_i^{-1}$.  
             For any Kripke Structure 
             $  {\cal K} =(S;X; J; 
             ( {\cal P}^j \ | \ j\in J) ;  \psi )$
and $i\geq 0$ we define $\phi ^{{\cal K}}_i : S \rightarrow \Omega_i(X,J)$     
 by $\phi_i^{{\cal K}} (s):= \phi^{{\cal K}} (s) \cap {\cal L}_i (X,J)= \pi_i 
  (\phi^{{\cal K}}(s))$.

             \vskip.2cm 

\noi 
             For every $0\leq i<\infty$ we 
consider the Kripke Structure  
         $\Omega_i= (\Omega_i; X;   J;(\overline{\cal F}_i^j\ |\ j\in J) ;
\overline{\psi}_i)$,  
          where $\overline {\psi}_i = 
          \overline {\psi} \circ \pi_i ^{-1} $ and for $i>0$  
           the partition $\overline {\cal F}^j_i$ of $\Omega_i$   
  is induced by the inverse images of  
 the function
  $\beta_i ^j: \Omega_i \rightarrow 2^{{\cal L}_{i-1}(X,J)}$ defined
  by $\beta_i ^j (w):= \{ f\in {\cal L}_{i-1}(X,J)\ | \ 
  k_j(f) \in w\}. $
  We  define $\overline {\cal F}^j_0=\{ \Omega_0\}$ for every 
   $j\in J$. 
\vskip.2cm 
   
           \noi  Now we consider $\Omega_i$ again as a canonical projective 
             image. 
            ${\cal G}_i$ is defined  to be the partition of $\Omega$ 
            induced by the inverse images of $\pi_i$,  
            ${\cal G}_i:= \{ \pi_i^{-1} (w)\ |\ w\in \Omega _i\}.$
             By the definition of $\Omega$, the join partition 
             $\vee_{i=1}^{\infty}  
             {\cal G}_i$ is the discrete partition of  
              $\Omega$, meaning that it consists of singletons.  
                      Let ${\cal F}^j_i $ be the  
                partition on $\Omega$, coarser than ${\cal G}_i$, 
                 defined by ${\cal F}_i^j:= 
                 \{ \pi_i^{-1} (B)\ |$    $\ B\in \overline 
                  {\cal F}^j_i\} .$
                From the definitions of the $\Omega _i$ and the 
${\cal F}^j_i$ it follows that 
                $\vee _{i=0}^{\infty} {\cal F}^j_i = {\cal Q}^j$. 

          \vskip.2cm

\noi 
Since 
$X$ and $J$ are  finite, 
there are several important properties of the Kripke Structure 
$\Omega_i$, all of which 
 are used in  this paper.
\vskip.2cm 
\noi 
(i) ${\Omega}_i$ is finite 
 for every $0\leq i < \infty$. (For a more general 
 statement, see Lismont and Mongin 1995.)  
          \vskip.2cm

\noi 
 (ii)  For every  $ w\in \Omega_i$ we can define  a formula 
       ${\bf f}( w)$  
       of depth  $i$ or less  
       such that $\alpha^ {\Omega_i}({\bf f}(w))=\{ w\} $, meaning that 
        the formula ${\bf f}(w)$ is true with respect to $\Omega_i$ only at 
         $w\in \Omega_i$.
 This follows from the finiteness  
         of $\Omega_i$. 
        For any subset $A\subseteq \Omega_i$   
         define ${\bf f}(A):= \vee _{w\in A}{\bf f}(w)$, a formula that is true  
          with respect to $\Omega_i$ only in the subset $A$.
          \vskip.2cm
\noi 
(iii)  It is easy to extend a member of  $\Omega _i$ to a member  
 of $\Omega_{i+1}$. 
Fix $0\leq i<\infty$ and $w\in \Omega_i$. For every  
 $j\in J$ define $\overline F_i^j$ by  
  $w\in \overline F_i^j\in \overline {\cal F}^j_i$.
  If  $(M_i^j\ | \ j\in J) $ are subsets 
  of $(\overline F_i^j \ | \ j\in J)$, respectively,    
    such that \newline 
   1) $w\in M_i^j$ for every $j\in J$, and\newline  
   2)  for every 
     $B\in {\cal G}_{i-1}\quad $   
   $\overline F_i^j \cap \pi_i(B)\not=\emptyset$ 
   implies that $M_i^j\cap \pi_i(B)\not= \emptyset$, \newline 
      then there is a unique  $v\in \Omega_{i+1}$ such that $\pi^{i+1}
 _i  (v)=w$  
       and for every $u\in \Omega_i\quad $ $\neg k_j \neg {\bf f}(u)\in v$ if and 
        only if $u\in M^j_i$. Furthermore, this is the only way 
 to extend a member of $\Omega_i$ to a member of $\Omega_{i+1}$; this is 
 Lemma 4.2 of Fagin, Halpern, and Vardi (1991).
        For any $i\geq 0$ and $v\in \Omega_k$ 
         with $k>i$ we  define $M_i^j(v) :=\{ u\in \Omega_i\ |\  
         \neg k_j \neg {\bf f}(u)\in v\} $. 
                  Notice  
              that if $w\in F\in \overline {\cal F}_i^j$ then 
             $M_{i-1}^j(w)$ is equal to $\pi^i _{i-1} (F)$, 
             which  could be a  
  proper subset of the member of $\overline {\cal F}_{i-1}^j$ that contains 
             $\pi^i _{i-1}(w)$.

          \vskip.2cm
 \noi 
 (iv) For every formula $f\in {\cal L}_i$   
and $l\geq i\quad \pi^{-1}_l( \alpha ^{\Omega _l} (f))
          = \alpha ^{\Omega}(f)  .$ This follows from (iii) and the 
 Completeness Theorem.   
        (See also Lemma 2.5 in Fagin, Halpern, and 
        Vardi 1991.)   
                                                                \vskip.2cm

\noi 
   (v) As a Kripke Structure, every $\Omega_i$ is connected.
   This was proven first by Fagin, Halpern, and Vardi (1991) and it can be 
 proven in several ways (for example from Proposition 1 of Simon, 1999). 
  
\subsection {The Common Knowledge of a Formula} 

Fagin, Halpern and Vardi (1991) investigated what happens when the agents 
 have common knowledge of a finite set of formulas, equivalently the 
 common knowledge of a single formula. 
Following  their definition  for 
 ``closed" and not wanting to create confusion with 
 topologically closed, for  all $i>0$  
we define a non-empty 
subset $A\subseteq \Omega _i$ 
to be {\em semantically closed} (Simon 2001)
 if for every $j\in J$,  every $B\in {\cal G}_{i-1}$ 
 and  every $w\in A$   if $\pi^{-1}_i(w)\subseteq F\in  {{\cal F}^j_i}$ and  
  $F\cap B\not=\emptyset$ then $F\cap B\cap \pi_i^{-1}(A)\not=\emptyset$. 
  Any non-empty subset of $\Omega_0$ is allowed to be  
  semantically closed.  Let $f\in {\cal L}$ be a formula with 
 $d=$ depth $(f)$. 
 Fagin, Halpern, and Vardi (1991) 
 proved that ${\bf Ck}(\{ f\})$ is not empty if and only 
 if the subset $\alpha ^{\Omega_d} (f)$ is a semantically closed subset 
 of $\Omega_d$ and that there exists a cell dense in ${\bf Ck}(\{ f\})$
 (equivalently $\underline {Ck}(\{f\})=F(C)$ for some cell $C$) if 
 and only if the Kripke Structure ${\cal V}(\alpha ^{\Omega_d}(f))$ is 
 connected.  
For all $i\geq d=$ depth $(f)$ we define 
 $\Omega^f_i:= \pi_i ({\bf Ck}(\{ f\} ))$; it follows from property 
 (iv) that $\Omega^f_i \subseteq \alpha^{\Omega_i} (E_{i-d}(f))$. 
Define $ \overline {\cal F}_i^j(f)$ by $\overline {\cal F}_i^j(f) 
:= \{ F \cap \Omega ^f_i
\ | \  F\in \overline {\cal F}_i^j\}$.
 Define the Kripke Structure $\Omega^f_i= (\Omega^f_i; J;  
(\overline {\cal F}_i^j(f)
 \ | \ j\in J); X ; \overline \psi_i|_f)$ 
 where 
  $ \overline \psi_i|_f$ refers to the restriction of $\overline \psi_i$ 
 to $\Omega_i^f$.
Likewise define ${\cal F}^j_i(f)$ by ${\cal F}^j_i(f) := 
 \{ \pi^{-1}_i(F)\ | \ F\in \overline {\cal F}^j_i(f)\} $ and 
 define ${\cal G}_i (f)$ by 
 ${\cal G}_i(f) := 
 \{ G\in {\cal G}_i\ | \ G\subseteq \alpha ^{\Omega} (E_{i-d}(f)).$ 
 (where $d$ is the depth of $f$). 

\vskip.2cm \noi 
Most importantly, Fagin, Halpern, and Vardi (1991) showed how to create  
 extensions of $\Omega^f_i$ to $\Omega ^f_{i+1}$ for all $i\geq d=$ depth 
 $(f)$, with only an additional requirement to the 
 rules of (iii) governing extensions from $\Omega _i$ to $\Omega _{i+1}$:
    We must require that the $\overline F^j_i$ are in $\overline 
 {\cal F}^j_i (f)$, which means that the 
  $(M_i^j\ | \ j\in J) $ are also subsets 
  of $ \Omega^f_i$. 
For the existence of such subsets is needed the semantically closed property. 
 The ability to extend establishes the equality  
    $\Omega^f_i = \alpha^{\Omega_i} (E^{i-d}(f))$ for all 
 $i\geq d=$ depth $(f)$.

%%%%%%%%%%%%%%%%%%%%%%%%%%%%%%%%%%%%%%%%%%     

\vskip.2cm 
\noi 
 Fix $w\in \Omega^f_i$ with $i\geq d=$ depth $(f)$ and 
with $\alpha^{\Omega_d} (f)$
 semantically closed, 
 and let $\overline F_i^j$ be the member of  $\overline 
 {\cal F}^j_i(f)$  
  containing $w$. The choice of 
$M_i^j(p_{i+1}(w))=\overline F^j_i$ for agent $j$ was 
  called the ``least-information" extension 
            in Fagin, Halpern, and Vardi (1991).
  Define $p^f_{i+1}(w)$ to be that unique  member of $\Omega^f_{i+1}$ 
   such that $\pi_i (p^f_{i+1}(w))=w$ and 
   $M^j_i(p_{i+1}(w))= \overline F_i^j$ 
   for every $j\in J$.   
  If $f$ is a tautology, then it was called the ``no-information" extension, 
 and in this case we write $p_{i+1}$ instead of $p^f_{i+1}$. 

\vskip.2cm 
\noi 
We define a formula $f\in {\cal L}(X,J)$ with $d=$ depth $(f)$ and 
 $|J|\geq 2$ to be {\em generative} if and only if $\alpha ^{\Omega_d}(f)$ is 
 semantically closed, ${\cal V} ( \alpha ^{\Omega_d}(f))$ is connected, 
 and 
 there exists more than one cell dense in ${\bf Ck} (\{ f\} )$, (meaning that 
 these cells are uncentered). 
In Simon (2001), Theorem 1 states that  the following are equivalent:\newline 
(a) the formula $f$ is generative, \newline 
(b)there is an uncentered 
 cell $C$ such that $F(C)= \underline {Ck}(\{ f\} )$, meaning that 
 there are  uncountably many such cells, (equivalently 
 uncountably many uncentered cells dense in ${\bf Ck}(\{
 f\})$),\newline  
(c) $\underline {Ck}(f)=F(C)$ for some cell, but $\underline {Ck}(f)$ 
  is not a maximal member of ${\cal T}$, \newline 
(d) there is a cell dense in ${\bf Ck}(\{ f\})$ that  is not finite.
\vskip.2cm 

\noi 
We will call any member of $\Omega_i$ an {\em atom}, or an {\em atom} 
 of $\Omega^f_i$ if 
 it also belongs to $\Omega^f_i$. 
%%%%%%%

            \section{Uncountably many cells
 with uncountable cardinality}
 
Our first  goal is to prove {\bf Theorem 1:} If  
 the formula $f$ is generative then there is a continuum of  
  distinct cells dense in ${\bf Ck}(\{ f\})$ such that there are infinitely 
 many possibility sets  with  continuum  
 cardinality. We prove Theorem 1 with something called the {\em alienated 
 extension.} The alienated extension is an alternation between different 
 ways to extend an element of the finite level 
 structure  $\Omega^f_i$, with long stretches of least information extensions 
  and long stretches of confirming the finite models $\Omega^f_k$ for 
 infinitely many $k$. 
\vskip.2cm 
\noi 
If $f$ is generative and $i\geq$ depth $(f)$ 
define an $ F\in \overline 
{\cal F}_i^j(f)$ to be {\em proto-generative} (for $f$)  if 
 there exists at least one $v\in \Omega^f_{i-1}$ such that the number 
 of members of $\Omega_i^f$  in  
$ F\cap \pi_i \circ \pi_{i-1}^{-1}  (v)$ is at least $2$; and 
 define such an $F\in \overline {\cal F}_i^j(f)$ to be {\em generative} 
 (for $f$) if for every  $v\in \Omega^f_{i-1}$ such that   
$ F\cap \pi_i \circ \pi_{i-1}^{-1} (v)\not= \emptyset$
  then 
 the cardinality of this intersection 
 is at least $2$.  Define an atom $w\in \Omega^f_i$ to 
 be  proto-generative (respectively generative) 
 for an agent $j$ if it is contained in a  member 
 of $\overline {\cal F}^j_i(f)$ that is proto-generative
 (respectively generative). 
  \vskip.2cm 
\noi 
 If $f$ is generative with depth 
 $d$ then there 
 must be a  proto-generative member of  $\overline 
{\cal F}^j_d(f)$ for some $j\in J$,
 since otherwise all extensions from $\Omega^f_d$ to 
 $\Omega^f_{d+1}$ would be determined uniquely, and the same would 
 be true  for all the following $\Omega^f_i$ for all $i>d$, and we 
 would have a contradiction to Theorem 1 of Simon (2001).  
\vskip.2cm 
\noi 
{\bf Lemma 1:} Let $f$ be generative and let $i\geq d=$ depth $(f)$. 
\vskip.2cm 
\noi 
{\bf (a)} If $F\in \overline {\cal F}^j_i(f)$
 is proto-generative then  every 
$G\in \overline {\cal F}^k_{i+1}(f)$ with $k\not=j$ 
 and  
 $\pi^{-1}_{i+1}(G) \cap \pi^{-1}_i (F)  \not= \emptyset$ is 
  also proto-generative. \vskip.2cm 
\noi 
{\bf (b)} Let $F\in \overline {\cal F}^j_i(f)$ be given. 
If every $G\in \overline {\cal F}_i^k(f)$ such that $k\not=j$ and 
 $G\cap F\not= \emptyset$ is  proto-generative, 
 then every $F'\in \overline {\cal F}^j_{i+1}(f)$ extending $F$ is 
 generative.   \vskip.2cm

\noi 
{\bf (c)}  If there are at least three agents then 
 there is a level $\hat i\geq d$ such that for all $k\geq \hat i$ 
 all $k$-atoms of $\Omega^f_k$  are generative for 
 all agents.  If there are exactly two 
 agents, then there is a level $\hat i\geq d$ such that 
 for all $k\geq \hat i$ any $k$-atom of $\Omega^f_k$ is  generative for 
 either one or the other agent.  \vskip.2cm 

\noi 
{\bf (d)} There is a level $\hat i$ such that 
   for all $ k \geq \hat i$ if  
 $F'\in \overline {\cal F}^j_{k+2}(f)$ was created from   
 the  use of the least information extension twice   
 from $F\in \overline {\cal F}_k^l(f)$  and $B$ is in $ {\cal G}_k$ with 
 $B\subseteq \pi_k^{-1} (F)$ then  
 the intersection $ F' \cap \pi _{k+2} (B)$ 
 has at least two elements. 
\vskip.2cm 
\noi {\bf (e)} If there are only two  agents and 
  $i$ is large enough so that all atoms of $\Omega^f_i$ are 
 generative for one or the other agent 
 but none of them are pro-generative for both agents then there must be one 
 agent $j$ such that all the atoms of $\Omega^f_i, \Omega^f_{i+2}, \dots$ 
 are generative  for $j$,  all the atoms of 
 $\Omega^f_{i+1}, \Omega^f_{i+3}, \dots$ are generative  for the other 
 agent  $j'\not= j$,  none of the atoms of   $\Omega^f_i, \Omega^f_{i+2}, \dots$ 
 are proto-generative for $j'$, and none of the atoms 
 of $\Omega^f_{i+1}, \Omega^f_{i+3}, \dots$ are proto-generative 
 for $j$. 
 \vskip.2cm 

\noi {\bf Proof:} 
\vskip.2cm 

\noi {\bf (a)} Let $F'\in \overline 
{\cal F}^j_{i+1}(f)$ be any extension of $F$ intersecting 
 $G\in \overline {\cal F}^k_{i+1}(f)$,
 and let $B\in {\cal G}_i(f)$ be any member such 
 that $\pi^{-1}_{i+1}(F')\cap \pi^{-1}_{i+1}(G)\cap B\not= \emptyset$. 
Since there is  at least two ways for Agent $j$ 
to extend $F$ that included the possibility of  $\pi_i (B)\in \Omega^f_i$ 
 (and because in extending  $\pi_i(B)$ the agents choose their sets $M^j_i$ 
 independently) 
we  conclude that $|G\cap \pi_{i+1}(B)|\geq 2$.  
\vskip.2cm 

\noi {\bf (b)} Let $B$ be any member of ${\cal G}_i(f)$ such that 
 $\pi^{-1}_{i+1}(F')\cap B \not= \emptyset$, and let $G\in 
\overline{\cal F}^k_i(f)$ 
 be such that $F\cap G$ contains $ \pi_i (B)$.  Because $G$ is proto-generative
 there must be at least two elements of $\Omega_{i+1}^f$ in 
 $F' \cap \pi_{i+1}(B)$. 
\vskip.2cm

 \noi 
{\bf  (c)}  If none of the atoms of $\Omega^f_d$ were proto-generative 
 then there would be only one way to extend all of these atoms to the next 
 level, and so on ad infinitum, and 
 that would contradict the assumption that $f$ is generative. Let 
 $v$ be a  atom of $\Omega^f_{d}$ that is not proto-generative 
  of adjacency distance one from 
 some proto-generative atom $w$ of $\Omega^f_d$, with the two atoms 
sharing membership of  $F\in \overline {\cal F}^j_d(f)$. Since 
 there is only one extension of $F$ to a member of $\overline 
 {\cal F}^j_{d+1}(f)$ it must hold that every extension 
of $v$ in $\Omega^f_{d+1}$  remains 
 adjacent to every extension of $w$ in $\Omega^f_{d+1}$. Because 
 $w$ is proto-generative, and so   for $j'\not=j$ the member of 
 ${\cal F}^{j'}_d(f)$ containing $v$ is proto-generative,
 it follows from Part (a) that every extension 
 of $w$ is also proto-generative. By induction on the adjacency 
 distance  it follows 
 that  if $w\in \Omega_d^f$ is of adjacency distance  $l$ from 
 a proto-generative $v\in \Omega_d^f$ then every extension 
 of $w$ in $\Omega_{d+l} ^f$ is also proto-generative. From the finiteness 
 of the adjacency diameter of $\Omega_d^f$  there 
 is an $l$ such that all atoms of  $\Omega_{d+l}^f$ are proto-generative. 
 Combining   Parts (a) and (b), 
 if $v\in \Omega^f_i$ is proto-generative  for agent $j$ then 
 any extension of $v$ in $\Omega^f _{i+1}$ is proto-generative for 
the other agent, any extension of $v$ in 
  $ \Omega^f_{i+2}$ is generative for agent $j$, 
any extension of $v$ in 
  $\Omega^f_{i+3}$ is generative for the other agent, and  so on.
  The claim concerning  more than two 
 agents is now transparent.   
\vskip.2cm  

\noi {\bf (d)} It follows directly  from  Parts (b) and (c). 
\vskip.2cm 
\noi {\bf (e)} That for the first level $\hat i$ 
  there are  generative atoms 
 for only one agent 
  follows from the connectedness of  $\Omega_{\hat i}^f$ and
 induction on the adjacency 
 distance.  
 The rest follows from   Parts (b) and (c).  

\hfill $\Box$
\vskip.2cm

\noi 
For any generative formula $f\in {\cal L}$ define 
 gen $(f)$ to be the first level $i\geq $ depth $(f)$ 
 such that  if there are two agents then 
 every  member  of  $\Omega ^f_i$ is  generative for 
 one or the other agent, and if there are at least three agents then 
 all members of $\Omega^f_i$ are generative for all agents.  
\vskip.2cm 

\noi 
For the rest of this section let a generative $f\in {\cal L}$ be fixed. 
Let $2_{\infty}^{\bf N_0}$ be the set of subsets of the whole  
numbers  ${\bf N_0}=\{0, 1,2, \cdots \} $ 
with infinite cardinality ($S\in 2_{\infty}^{\bf N_0}$ implies 
$S\subseteq {\bf N_0}$ and  
 $|S|=\infty$).      For every pair $i,k\in S$ with $k\geq i\geq $ depth 
 $(f)$ we will 
  define a map $p^{S,f}_k :\Omega^f_i \rightarrow \Omega^f_k$.
    If $i\in S\in 2_{\infty}^{\bf N_0}$ define $n_S(i):= \mbox { inf } 
  \{ k\in S\ | \  k>i \}$, the first member of $S$ strictly larger than $i$. 
  If $i\in S$ and $w\in \Omega^f_i$ 
then   define $p^{S,f}_{n_S(i)}(w) :=  \phi ^{\Omega^f_i}_{n_S(i)} (w)$   
  and define $p_i^{S,f}(w) :=w$.
       $p^{S,f}_{n_S(i)}(w)$ is an extension of $w$, 
        meaning that $\pi^{n_{S}(i) }_i
 (p^{S,f}_{n_S(i)}(w))=w$ (Lemma 1 of Simon 2001).
      For every $k\in S$ and  
       $w\in \Omega^f_i$ with  $k\geq i\in S$ and $p^{S,f}_k(w)\in \Omega^f_k$
already defined,  
      define $p^{S,f}_{n_S(k)}(w)$ 
       to be  $ p^{S,f}_{n_S(k) } ( p^{S,f}_k (w))$.   
    Lastly, for all $i\in S\in 2_{\infty}^{\bf N_0}$ and  
     $w\in {\Omega^f}_i$ define $p^{S,f}:\Omega^f_i \rightarrow \Omega^f$ by  
     $$p^{S,f} (w):=\bigcap_{l\in S, \ 
    l>i}\pi^{-1}_l \circ p_l^{S,f}  
     (w).$$  
      For any $i\in S\in 2_{\infty}^{\bf N_0}$ and 
           $w\in \Omega^f_i$ we call  
            $p^{S,f}(w)$ the  alienated extension 
            of $w$ with respect to S and $f$.     
 Define $p^S$ to be $p^{S,f}$ for any tautology $f\in {\cal L}$.            
\vskip.2cm 

\noi 
           An  alienated extension involves an     
             infinite number of least-information extensions. For all 
              $0\leq i<\infty$ and $w\in \Omega^f_i$ it is easy to confirm 
 that        $\phi_{i+1} ^{\Omega^f_i}(w)=p^f_{i+1}(w)$, meaning also 
                that $p^{{\bf N_0},f}$ is the infinite repetition of the 
                 least-information extension. 
  Define the map $p^f$ to be $p^{{\bf N_0},f}$ and $p$ to be $p^{\bf N_0}$.

\vskip.2cm 
\noi 
For any $S\in 2^{\bf N_0}_{\infty}$ and positive $k$ define $n_S^k(i)$ by 
 $n_S^1(i)=n_S(i)$ and $n_S^k(i)=n_S\circ n_S^{k-1} (i)$. 
\vskip.2cm

     \noi   {\bf Lemma 2:} If $S\in 2_{\infty}^{\bf N_0}$
 and $f$ is generative, then 
        all alienated extensions with respect to $S$  and $f$ share the same  
         dense  cell of ${\bf Ck}(\{ f\})$.\vskip.2cm  
   
  \noi 
 {\bf Proof:} If $i\geq $ depth $(f)$ and both 
     $w$ and $w'$ are  members of $\Omega^f_i$ such that both are contained  
      in the 
      same member  of $ \overline {\cal F}^j_i(f)$, then 
       from induction and the definition of $\phi^{\Omega^f_i}$ 
$p^{S,f} (w)$ and  
       $p^{S,f}(w')$  are both contained in  
       the same member of ${\cal Q}^j$, the limit of the 
$\overline {\cal F}_i^j$. 
\vskip.2cm 

\noi   Now, given any  $i,k\in S $ and $b\in \Omega^f_i$ and 
   $d\in \Omega^f_k$, 
   the adjacency distance between $p^{S,f}(b)$ 
    and $p^{S,f}(d)$ in $\Omega^f$ is no more than the adjacency distance 
     between $p^{S,f}_{\max (i,k)}(b)$ and 
     $p^{S,f}_{\max(i,k)}(d)$ in $\Omega^f_{\max (i,k)}$. 
  The rest follows by the connectedness  
  of $\Omega^f_i$ for every $i\geq $ depth $(f)$.
  \hfill $\Box$ \vskip.2cm

 \vskip.2cm 
\noi 
     Given any generative formula $f$ define the formula 
     $g^f_i:=
{\bf f}(\phi_{i+1}^{\Omega^f_i}(\Omega^f_i))$ of depth $i+1$, the formula 
      true in  $\Omega^f_{i+1}$ only in the image 
      $\phi_{i+1}^{\Omega^f_i}(\Omega^f_i)$.
       As we will see, the common knowledge of $g^f_i$ is closely linked 
       to the Kripke Structure $\Omega^f_i$ 
       (see also Theorem 4.23 of Fagin, Halpern, 
        and Vardi 1991).
      \vskip.2cm

    \noi {\bf Lemma 3:} The formula  $g^f_i$ is common knowledge 
        in the Kripke Structure $\Omega^f_i$. 
    If  $i\in S\in 2_{\infty}^{\bf N _0}$, $i\geq$ depth $(f)$,  
    and $i+1, i+2, \dots ,i+l \not\in S$ for some $l\geq 1$,  
     then $p^{S,f} (\Omega^f_i)  
     \subseteq \alpha ^{\Omega^f}  
      (E^l(g^f_i))$, and the same holds for
  $\phi_{i+l+1} ^{\Omega^f_i}$ applied to any element 
 of $\Omega^f_i$. 
\vskip.2cm

     \noi {\bf Proof:}  
     Because $\Omega ^f_i$ is finite and connected,
 $\phi^{\Omega^f_i}(\Omega^f_i)$ is a cell. Because 
     $\phi^{\Omega^f_i}(\Omega^f_i)\subseteq 
      \alpha^{\Omega^f} (g^f_i)$,  Property (iv) and  
            Lemma A  imply that $g^f_i$ 
            is  common knowledge in 
     the cell     $\phi^{\Omega^f_i}(\Omega^f_i)$.
       If $E^l(g^f_i)$, a formula of depth $i+l+1$, 
      were  not true 
        at some point of $\phi^{\Omega^f_i}_{i+l+1}(\Omega^f_i)$ then also   
         by  Property (iv) we must have that 
         $g^f_i$ is not common knowledge at some point of 
                 $\phi^{\Omega^f_i}(\Omega^f_i)$, a contradiction.
     \vskip.2cm 

\noi   
       By Simon (2001) $\Omega^f_i$ and $\phi^{\Omega^f_i} (\Omega ^f_i)$ 
 are equivalent as Kripke Structures, so that 
 $g^f_i$ is also common knowledge 
        in the Kripke structure  $\Omega^f_i$. 
     \hfill $\Box$ 
         \vskip.2cm

\noi {\bf Lemma 4:} If $f$ is generative and   $i\geq $ gen $(f)$
 then  $ E g^f_i$ is not true at 
 any extension   of $p^{f}_{i+2}(\Omega^f_{i+1})$.       \vskip.2cm 

\noi 
{\bf Proof:}    Because $\phi_{i+2}^{\Omega^f_{i+1}}(w)= p^f_{i+2} (w)$ for 
 any $w\in \Omega^f_{i+1}$, 
given  
$w\in   \Omega^f_{i+1}$ it suffices to find 
 some $j\in J $ such  that   
$k_j ( g^f_i)$ is not true at 
 $p^f_{i+2}(w)$. Let $j\in J$ be such that 
  $w\in  \Omega^f_{i+1}$
 is generative for Agent $j$ and 
let $w$ be in $ F\in \overline {\cal F}^j_{i+1}(f)$.
 For every $v\in \Omega^f_{i}$ there is only one member of $\Omega^f_{i+1}$ in 
    the subset $\pi_{i+1}(\pi^{-1}_i(v))$
  where $g^f_i$ is true.
    But  for 
     $v:= \pi^{i+1} _i  (w)\in \Omega_i^f$ 
  there are at least two 
   $u\in \Omega^f _{i+1}$ with  
  $u\in F \cap \pi_{i+1}(\pi_{i}^{-1}(v))$ (including 
 at least the possibility of $u=w$).  
 The $F'\in \overline {\cal F}^j_{i+2}(f)$ containing 
 $p^f_{i+2}(w)$ must have a 
  non-empty intersection with  $\pi_{i-2} \circ \pi_{i+1}^{-1}(u)$ for all 
  $u\in F \cap \Omega^f_{i+1} $, 
  and therefore $F'$ is   
   not contained in $\alpha^{\Omega^f_{i+2} }( g^f_i)$. 
                                                    \hfill $\Box$\vskip.2cm

 \noi 
 {\bf Proof of Theorem 1:} 
 Define a map $\beta :2^{\bf N_0}\rightarrow 2 
 _{\infty}^{\bf N_0}$ by \newline  $\beta (S) :=  
 \{0, 1,2,4,8,\cdots \} \cup \{ 2^i+1, \cdots, 2^{i+1}-1\ | \  i\in 
  S\} $. \newline   
  Define an equivalence relation  
  on $2^{\bf N_0}$ by $S\sim T$ if and only if there exists an 
  $m\in {\bf N_0}$ such that $S\backslash 
  \{ 0, 1,2,\cdots m\} =T\backslash \{0, 1,2,\cdots, m\}$. 
   The co-sets of this equivalence relation have 
   the cardinality of the continuum. 
\vskip.2cm 

\noi 
Let $k$ be such that  $2^k \geq $ gen $(f)$ and let $d=2^k$. 
 Due to Lemma 2, it suffices to show for some $w\in \Omega^f_d$ that if 
   $S$ and $T$ are both subsets of ${\bf N_0}$ with $S\not\sim T$ 
   then $p^{\beta (S),f} (w) $ does not 
    share the same  cell 
     as $p^{\beta (T),f}(w)$. 
      For the sake of contradiction, let us suppose  
      that the adjacency distance in ${\bf Ck}(\{ f\})$ between 
      $p^{\beta (S),f}(w)$ and $p^{\beta (T),f}(w)$ equals a finite number 
       $l<\infty$.  Because 
      $S\not\sim T$ there exists an $i>\max (log_2((l+2)), k)$
 such that $i\in S$  
       and $i\not\in T$, or vice versa.  By symmetry, 
        let us assume that $i\in S$ and $i\not\in T$.  
           Lemma 3 applied to $p^{\beta (T),f} (w)$  
          implies that $p^{\beta (T),f}(w)\in \alpha^{\Omega^f} 
       (E^{l+1} g^f_{2^i})$. But because the adjacency-distance   
      between  $ p^{\beta (S),f}(w)$ and $p^{\beta (T),f}(w)$ is $l$ we have   
         that $p^{\beta (S),f} (w) \in \alpha^{\Omega^f} (E(g^f_{2^i}))$, 
          a contradiction to Lemma 4.\vskip.2cm 
\noi If $S$ is infinite (all but one of the  uncountably 
 many equivalent relations) then 
   any    possibility sets containing an alienated extension  
 is  homeomorphic to a Cantor 
 set. This  follows from Lemma 1d 
  and the fact that for every such $S$  there are arbitrarily long strings 
of the least information extension applied.  
 \hfill $\Box$
\vskip.2cm 
 \noi Because there are only countably many alienated 
  extensions $p^{\beta (S),f} (w)$ 
 for any fixed $S$ (because there are only countably many atoms $w$ 
 on any level) we haven't shown that there are uncountably many cells 
 where   
  {\bf all}  the possibility sets are equivalent 
 to Cantor sets. That would be a more difficult claim to prove and would 
require a better understanding of what happens
 at any finite adjacency distance 
 from a point created by an 
  alienated extension. We leave it as an open question: 
\vskip.2cm 
\noi {\bf Question 5:} For any generative 
  formula $f$  does there exist uncountably many 
 cells  dense in ${\bf Ck}(\{ f\})$
 where all the possibility sets are equivalent 
 to Cantor sets? 
 \vskip.2cm 
 \noi {\bf Corollary 1:} If $f$ is generative and $C$ is a cell 
 of finite fanout where $f$ generates all 
 the formulas held in common knowledge 
   then there cannot be a
 uniform bound on the size of the possibility sets in $C$. 
 \vskip.2cm 
\noi Suppose for the sake of contradiction that $n$ is a uniform bound 
 on the size of the possibility sets in $C$. 
Let $F$ be any possibility set topologically equivalent to a Cantor set 
 in a cell $C'$ that is also dense in  ${\bf Ck}(\{ f\})$ and let 
 $z$ be any point of $F$.  Let $j$ be the agent such that 
  $F$ is a possibility set for $j$. 
 Because $F$ is infinite, there will be $n+1$  
 mutually exclusive  formulas $g_0,  \dots , g_n$ such that 
 $\neg k_j \neg g_i$ is true at $z$ for every choice of $i=0, \dots n$ ($\neg
 k_j \neg g_i$ meaning that $j$ considers the validity
 of $g_i$ to be possible). Because $C$ is also dense in  ${\bf Ck}(\{ f\})$ 
 there will be a sequence $z_1, z_2, \dots$ of points in $C$ converging 
 to $z$. At some level $L$ all the  $\neg k_j \neg g_i$ will be valid 
 at all $z_l$ for all $l\geq L$ and $0\leq i \leq n$. Because the $g_i$ 
 are mutually exclusive the possibility set for agent $j$ 
  containing any such point $z_l$ will have cardinality at least 
 $n+1$, a contradiction.  
   
\section {Finite Fanout}

Now we construct uncountably many cells of finite fanout dense 
 in ${\bf Ck}(\{ f\} )$ for any generative $f\in {\cal L}$. Let 
 such an $f$ be fixed for the rest of this section. 
         \vskip.2cm 
\noi The proof of Theorem 2 is quite different from that of Theorem 1. 
 Again we alternate how extensions are performed, but much more toward 
 confirming finite models and also selectively within any stage, so that 
 some points of $\Omega_i^f$ can be extended by way of the least information 
 extension and other points of $\Omega_i^f$ according to a previous 
finite model.
\vskip.2cm 

 \noi 
For any $w\in \Omega^f_i$ define $F^j_i(w)$ to be that member of 
 $\overline {\cal F}^j_i(f)$ containing $w$.
\vskip.2cm 

   \vskip.2cm 
\noi 
From Proposition 2 of Simon (1999), with finitely many agents 
 a cell is compact if and only if 
 it has finite diameter.  Therefore by Theorem 1 of  Simon (2001) 
 all cells dense in 
${\bf Ck}(\{ f\})$ do not 
 have finite diameter, and this implies also that there is no bound on 
 the diameters of the $\Omega^f_i$. 
\vskip.2cm 
\noi 
   Let $S\in 2_{\infty}^ {\bf N_0}$ satisfy \newline  
   1)  $\inf S> $ gen $(f)+8$,\newline 
   2)     for every $i\in S$   the differences 
 $n_S(i)-i$ (the next member of $S$ after $i$ minus $i$)
 start at  at least 5 and  
 are  strictly  increasing,
   \newline 3) for every $i\in S$ 
 the adjacency diameter of $\Omega^f_{n_S(i)}$
 is  strictly greater than twice the size of the set 
$\{ k \in S\ | \ k\leq i\}$ plus $3$, and 
\newline 4)  $2^{(n_s(i)-i-1)/2} $ is strictly greater 
 than  the cardinality of  $\Omega^f_i$.  
\vskip.2cm 

\noi 
Let $T$ be any infinite subset of $S$ with $\inf T =\inf S$. 
   For every such $T$ and 
 $i\geq \inf T$ we define inductively two subsets $A_i$ and 
    $B_i$ of $\Omega^f_i$.  
 If there are only two agents and  the levels beyond gen $(f)+2$ alternate 
 between being all generative for one agent and not proto-generative 
 for the other agent (Lemma 1e) then define  $w_0$  to be any
 element of $p^f_{\inf T} (\Omega^f_ {\mbox {gen} (f)+2})$.
 If there is an atom $v$ 
 of $\Omega^f_ {\mbox {gen} (f)+2}$ 
  that is generative for both agents then let $w_0$  be  $p^f_{inf T} (v)$, in 
either case 
 $w_0$ is an  application of 
     the least information extension   
  from  level   $\mbox {gen} (f)+2$ to the level $\inf T$. 
   Define  $B_{\inf T}$ to be the singleton $\{ w_{0}\}$
 and $A_{\inf T}=\emptyset$. 
     We assume that 
      $A_k$ and $B_k$ have been defined for all ${\inf T}\leq k<i$, and 
         show how to define $A_i$ and $B_i$. 
       First, we  define a extension function $\gamma_{i} :A_{i-1} 
       \cup B_{i-1} 
        \rightarrow A_i \cup B_i$ for all $i>\inf T$; it suffices 
         to determine the sets $M^j_{i-1}(\gamma_i(w))$.   
                    If $w\in A_{i-1}$ 
            then $M_{i-1}^j(\gamma_i (w)):= (A_{i-1}\cup B_{i-1}) \cap 
             F^j_{i-1}(w)$. If $w\in B_{i-1}$  and  
             $F^j_{i-1}(w)$ contains
                some member 
               of $A_{i-1}$, then 
             $M_{i-1}^j(\gamma_i (w)):= (A_{i-1}\cup B_{i-1}) \cap 
             F^j_{i-1}(w)$; otherwise   
             if $A_{i-1}\cap F^j_{i-1}(w)=\emptyset$, 
               then 
             $M_{i-1}^j(\gamma_i (w)):= F^j_{i-1}(w)$. 
       If $i\in T$ we define $B_i$ to  
       be the set 
        $\{ p^f_{i} (w) \ |\ w\in \Omega^f _{i-1}\backslash 
        (A_{i-1}\cup B_{i-1})$,  $w\in F^j_{i-1}(b)$
 for some $b\in B_{i-1}$ 
         and $j\in J$ with   
          $F^j_{i-1}(b)\cap A_{i-1}=\emptyset\} $     and 
          we define $A_i$ to be 
           the set $\gamma_i (A_{i-1}\cup B_{i-1})$. 
         If $i\not\in T$ we define $B_i$ to be the set 
          $ \gamma_i (B_{i-1}).$ 
             and we define $A_i$ to be the set 
                           $\gamma_i (A_{i-1})$. 
         For any $i> \inf T$, $l\geq 0$, and $w\in A_{i-1}\cup B_{i-1}$ 
         we define 
         $\gamma_{i+l}(w) = \gamma_{i+l} \circ \dots  
          \circ \gamma _i (w)$ and we define 
          $\gamma (w) := \cap_{k=i}^{\infty} \pi^{-1}_k\gamma_k(w)$. 
          We define $C$ to be $\{ \gamma (w)\ |\ i\geq \inf T, w\in A_i\}$. 
\vskip.2cm
\noi  If there are only two agents, notice from Lemma 1 
  that  any 
 atom of $\Omega^f_i$ with $i\geq \inf S$ and  within an adjacency distance 
 of  $3$ from $B_i$  is either generative for both 
 agents or it is generative for one agent and not proto-generative 
 for the other agent -- and the same holds for all extensions 
 of this atom to higher levels. It would be  nice, if possible, to prove 
 that if  any atom of $\Omega^f_{\mbox {gen}(f)+2}$ is generative 
 for both agents then there is a level for which all atoms are generative 
 for both agents. Such an argument was easy for non proto-generation for 
 both agents, since  the extensions of such atoms 
  were determined and  not 
 ``running away''. However  proto-generation with one but not the other agent 
 still allows for distinct extensions of the same atom to gain  distance
 from each other. The question is whether they can do so fast enough 
 to avoid mutual generation.\vskip.2cm

\noi  
    {\bf Lemma 5 :} The extension function 
 $\gamma_i$ is well defined for every $i\geq \inf T$     
 and if   $b\in B_i$ is adjacent in $\Omega^f_i$ to $a\in A_i$,  
  sharing the same member of $\overline {\cal F}^j_i(f)$, and $k$ is the  
   largest member of $T$ less than or equal to $i$, then 
 $a=\gamma_i (b')$ for some $b'\in B_{k-1}$  
  with  $F^j_{k-1}(b') \cap A_{k-1}=\emptyset$. 
\vskip.2cm        

\noi 
{\bf Proof:} We prove both claims together by induction on $i$.
       $\gamma_{\inf T +1} (w_0)=p^f_{\inf T+1} (w_0)$ is well defined. 
                We assume that $\gamma_k$ is 
well defined for all $\inf T+1 \leq 
                k<i$.
 Let  $w\in A_{i-1}\cup B_{i-1}$, and 
for any given $j\in J$ let us assume that $v\in \Omega^f_{i-2}$ 
satisfies $\pi_{i-2}^{-1} (v) \cap \pi_{i-1}^{-1}
 (F^j_{i-1}(w)) \not= \emptyset$. 
                              We need
 to show that $\pi_{i-1}(\pi^{-1}_{i-2}(v)) 
                              \cap M^j_{i-1}(\gamma_i (w))\not= \emptyset$. 
          If $i-1\not\in T$ and $i>\inf T$ 
          then the   well definition of $\gamma_{i-1}$  shows the same  
           for $\gamma_i$, so  
 for the following cases, 
     we assume that $i-1\in T$.
\vskip.2cm 
\noi 
{\bf Case 1; $w\in A_{i-1}$ and $v\in A_{i-2}\cup B_{i-2}$:}  
$\gamma_{i-1}(v)$ is in $ F^j_{i-1}(w)$ 
because $v$ and $\pi^{i-1} _{i-2} (w)$ share   
  the same member of $\overline {\cal F}^j_{i-2}(f)$. 
\vskip.2cm 
\noi 
 {\bf Case 2; $w\in A_{i-1}$ and      $v\not\in A_{i-2}\cup B_{i-2}$:} 
       This is possible only if  
        $\pi^{i-1} _{i-2}(w) \in     
        B_{i-2}$ and $F^j_{i-2}(\pi^{i-1} _{i-2}(w)) \cap A_{i-2} 
        = \emptyset$. Since $v\in 
 F^j_{i-2}(\pi^{i-1} _{i-2}(w))$ 
         we have $p^f_{i-1}(v)\in B_{i-1}\cap F^j_{i-1}(w)$. 
\vskip.2cm 
\noi 
         {\bf Case 3; $w\in B_{i-1}$ and $ F^j_{i-1}(w) \cap A_{i-1}\not= 
         \emptyset$:}  Let $a\in      F^j_{i-1}(w) \cap A_{i-1}$. 
          By the second part of the induction hypothesis 
          $\pi_{i-2}^{i-1}  (a) \in B_{i-2}$ 
           with $F^j_{i-2}(\pi_{i-2}^ {i-1}(a))\cap A_{i-2}=\emptyset$.  
           It follows that $v$ is in
 $ F^j_{i-2}(\pi_{i-2}^{i-1}(a))$ and  whether or not 
 $v$ was in $B_i$ that there is an extension of $v$ in $A_{i-1} \cup B_{i-1}$.
\vskip.2cm 

\noi 
{\bf Case 4; $w\in B_{i-1}$ and $F^j_{i-1}(w) \cap A_{i-1}=\emptyset$:} 
Since $w\in p^f_{i-1} (\Omega^f_{i-2})$ we have that 
$p^f_{i-1} (v) \in F^j_{i-1}(w)$.  
\vskip.2cm 

\noi 
              For the second part of the claim, suppose for the sake of 
               contradiction that 
$b':=\pi_{k-1}^i  (a) \in A_{k-1}$. 
   $b'$ shares the same member of  
    $\overline {\cal F}^j_{k-1}(f)$ with $c:= \pi_{k-1}^i  (b) 
    \in \Omega^f_{k-1} \backslash 
   (A_{k-1}\cup B_{k-1})$. 
      For every $j\in J$ and $D\in {\cal G}_{k-2}(f)$  if 
$\pi^{-1}_k (F^j_k (\gamma_k (b')))$  
     intersects $D$ then it intersects $D$ in only 
      one member of ${\cal G}_{k-1}(f)$.  If $b'$ is generative for 
 $j$ then by Lemma 1 it  
      is different from  
      $\pi^{-1}_k(F^j_k (p_k^f (c)))=
\pi^{-1}_k(F^j_k (\pi_k^i i (b)))$, a contradiction.
   If there are only 
 two agents and  $b'$ is  
not generative for $j$ then by the well definition of $\gamma_{k-1}$ 
 $c$ must have been in $A_{k-1}\cup B_{k-1}$, also a contradiction. 
 So we conclude that $b'$ was in $B_{k-1}$.
      Furthermore, if $F^j_{k-1}(b') \cap A_{k-1} \not= \emptyset$ then 
        either $\pi_k^i   (a)$
 and $\pi_k^i   (b)$ would 
        not share the same member of $\overline {\cal F}^j_k(f)$ (the 
 case of $F^j_{k-1}(b')$ generative) or from the well 
 definition of $\gamma_k$ $\pi_{k-1}^i  (b)$ 
 would be in $A_{k-1}\cup B_{k-1}$ (the 
 case of $F^j_{k-1}(b')$ not generative), both contradictions.
      \hfill $\Box$
      \vskip.2cm 

\noi  The second part of Lemma 5 shows that 
 if $i\in T$, $b\in B_{i-1}$ and $\gamma_i (b)=a\in A_i$ then for 
 every $k\geq 1$ the only members of $A_{n^k_T(i)}\cup B_{n^k_T(i)}$ 
 adjacent to $\gamma_{n^k_T(i)}(a)$ 
  are already in the set $\gamma _{n^k_T(i)}(A_i \cup B_i)\subseteq 
 A_{n^k_T(i)}$. 
  Therefore  $C=\{ \gamma (w)\ |\ i\geq \inf T, w\in A_i\}$ 
  is a cell with finite fanout.
                \vskip.2cm
              
\noi 
{\bf Lemma 6:}
 If   $w$ and $w'$ are  in $\Omega^f_i$ for some  $i\geq \inf T$ but neither  
are  in     
   $B_i$ then 
 there is an adjacency path 
  $w=w_0 , w_1, \dots , w_q=w'$ in $\Omega_i^f$ 
 between $w$ and $w'$  such that 
    $w_m\not\in B_i$ for all $1\leq m\leq q-1$. 
                \vskip.2cm 

\noi 
{\bf Proof:} We proceed by induction on $i$.   Consider an adjacency path 
 $v_1, v_2, \dots, v_l$ in $\Omega^f_{i-1}$ 
 with  
  $v_1=\pi_{i-1}^i (w)$ and 
  $v_l=\pi_{i-1}^i (w')$.  We assume that 
  $v_k$ and $v_{k+1}$ share the same member $F^{j_k}_{i-1}$ of 
   $\overline {\cal F}^{j_k}_{i-1}(f)$ for all $1\leq k\leq  l-1$, and 
 that $j_k\not= j_{k+1}$ for every consecutive pair $k,k+1$.  
   We will define an extension $w_k \in \Omega^f_i$ of $v_k$  for every 
     $k$ such that 
      $M^{j_k}_{i-1} (w_k) = M^{j_k}_{i-1} (w_{k+1}) $ 
      is a subset  
$F_k^{j_k}\in \overline {\cal F}^{j_k}_{i-1}(f)$
 containing both  $v_k$ and $v_{k+1}$ 
       with a non-empty intersection with  $\pi_{i-1}(D)
       $ for every $D\in {\cal G}_{i-2}(f)$ 
        with $\pi^{-1}_{i-1}(F^{j_k}_k) 
        \cap D\not=\emptyset$ and  
                      $M^j_{i-1}(w_1) = 
            M^j_{i-1}(w)$ for at least one $j\not= j_1$ and 
           $M^j_{i-1}(w_l) = 
            M^j_{i-1}(w')$ for at least one $j\not= j_{l-1}$. 
            The path 
            $w, w_1,  
             \dots , w_l, w'$ will be 
  an adjacency path connecting 
              $w$ and $w'$, allowing  possibly for the identity of $w$ and 
               $w_1$ or of $w'$ and $w_l.$ We must 
show that these extensions can be done so that for all $1< k < l$ 
 no extension  $w_k$ is in $B_i$. 
\vskip.2cm 

 \noi 
{\bf Case 1; there are at least three agents:}   
Looking at any $v_k$, let  $j\in J$ be any agent other than $j_k$ or 
 $j_{k-1}$. Since  all levels involved are generative for all 
 agents and the choice 
 of $M^j_{i-1} (w_k)$ doesn't affect the connectivity,
 the selection of $M^j_{i-1} (w_k)$ can be made 
 so that $w_k$ is not a least information  extension and therefore
 not in   
 $B_i$. \vskip.2cm 

\noi 
{\bf Case 2; there are only two  agents ($J=\{ 1,2\}$) and 
 $i\not\in T$:} By $i\not\in T$ the only members of $B_i$ are 
 extensions of members of $B_{i-1}$. Since all members of $B_{i-1}$ are 
 adjacent to members of $A_{i-1}$, by the induction assumption we can 
 assume that at most $v_1$ and $v_l$ are in $B_{i-1}$. So if 
  $v_1$ and $v_l$ are not in $B_{i-1}$, we are done by induction 
 on the stages. 
  Without loss of generality  assume for now  that 
  $v_1=\pi_{i-1}^i   (w)$ is in $B_{i-1}$.
Let $m$ be the largest member of $T$ that is less than $i$. If $m$ is greater 
 than $\inf T=\inf S$ 
we know from Lemma 5 that there is a $\hat j\in \{ 1,2\}$ such  
that $v_1$ is a member of $B_{i-1}$ and  $\pi_{m-1} ^{i-1}  
(v_1)$ shared the same member of 
 $\overline {\cal F}^{\hat j}_{m-1}(f)$ with a $b\in B_{m-1}$ but with  no 
 member of $A_{m-1}$. In this case define $b'= \gamma_{i-2}(b)$
  (with $b'=b$ if $m=i-1$). Otherwise let $b'$ be any other member  of 
 $\Omega_{i-2}^f$ that shares a member of $\overline {\cal F}^{\hat j}_{i-2}(f)$ 
 with $\pi_{i-2}^i (v)$ (for which there must be several since 
 $w_0$ was created from several applications of the non-information extension).
 Because $\pi^{-1}_{i-1}(F^{\hat j}_{i-1}(v_1))\cap 
 \pi^{-1}_{i-2}(b')\not=\emptyset$ it follows that 
$F^{\hat j}_{i}(w)$ contains some extension $u\in \Omega^f_i$ 
 of this $b'$, and furthermore no member of either $B_i$ or $B_{i-1}$
 is an extension of this 
 $b'$.  Therefore, since 
 $w$ and $u$ are adjacent, we can replace $w$ by
  $u\not\in B_i$, do the same for 
 $w'$  if necessary, and repeat with the induction assumption 
 with the added  assumption 
 that  none of the $v_k$ are in $B_{i-1}$.  
     \vskip.2cm

\noi 
{\bf Case 3; there are only two agents ($J=\{ 1,2\}$) and $i\in T$:}
All members of $B_i$ are  created as extensions 
 of points not in $B_{i-1}\cup A_{i-1}$ and  they are all in 
 $p^f_i (\Omega^f_{i-1})$. 
  If possible, for every $k$  
     let $M^{j_k}_{i-1} (w_k) = M^{j_k}_{i-1} (w_{k+1}) $ 
      be any proper subset of $F_k^{j_k}(v_k)$
 containing both $v_k$ and $v_{k+1}$, meaning 
 that if this is possible then we
  keep both $w_k$ and $w_{k+1}$ out of $B_i$. 
         On the other hand, if  
         $M^{j}_{i-1} (w_k) =   
         F^{j}_k(v_k)$ is  forced 
 for both agents $j$ then we will show that the 
 so defined $w_k$ is also not   in  $B_i$. 
                     \vskip.2cm  
     
\noi    
{\bf Case 3a; 
 and $1<k<l$:} 
Either $v_k$ is generative for Agent $1$ or Agent $2$. 
Without loss of generality assume that $v_k$ is generative for 
 Agent $1$, with $F = F_{i-1}^1 (v_k) \in \overline {\cal F}^1_{i-1}  (f)$ 
connecting 
 $v_k$ to $v^*$, equal to either $v_{k-1}$ or $v_{k+1}$.  
 Because $F$ is generative, if there were only  way  
 to define $M_{i-1}^1(v_k)$ so as to  include both $v_k$ and $v^*$ 
 it must follow that     $F^1_{i-1}(v_k)=\{ v_k, v^*\}$ with 
 both $v_k$ and $v^*$ extensions of the same  $u\in \Omega^f_{i-2}$ 
 (with $v^*= v_{k+1}$ if $j_k =1$ 
 and $v^*=v_{k-1}$ if $j_k =2$). Looking at  
 what is necessary for $w_k$ to be in $B_i$, since $F$ projected to 
 the $i-2$ level contains only one element and either the previous 
 number  in $T$ is below the $i-4$ level or $i=\inf T$ and $w_0$ was created 
through several applications of the least information extension, 
it is also necessary for 
 there to be some  $u'\in  F^2_{i-2} (u)$ other than $u$ (with the 
 possible choice of  $u'\in B_{i-2}$ 
 if $i> \inf T$; there will be many more than one other member of 
 $F^2_{i-2} (u)$, but one other suffices for the argument). 
 If $v_k$ were generative for Agent 2 we 
 would be able to avoid $w_k$ in $B_i$ from the
 choice of  $M_{i-1} ^2(v_k)$.  So we have to assume that $v_k$ is 
 not generative for Agent 2 and therefore  from Lemma 1  that  
 $F^2_{i-2} (u)$ is generative. 
Let $H$ be $\pi_{i-3}^{-1} \circ \pi^{i-2}_{i-3} (u')$.   
  There must be at least three members of $F$ corresponding 
 to the three non-empty subsets of  $\pi_{i-2}(H) \cap \Omega^f_{i-2}$ 
 combined with $u$, a contradiction to our assumption that 
 $F=F^1_{i-1}(v_k)$ contained only two elements. Therefore 
 there was more than  
 one way to connect our extensions of $v_k$ and  $v^*$, one of those ways  
 avoiding membership of $w_k$ in $B_i$.     
 \vskip.2cm 

     \noi    {\bf Case 3b;  $l>1$ and 
  $k=1$ or  $k=l$):} Without loss of generality we assume 
 that $F\in \overline {\cal F}_{i-1}^1 (f)$ contains 
 both  $v_1$  and $v_2$. 
 As with Case 3a, if $v_1$ were generative for Agent $1$  there would be 
 no alternative to  $M_{i-1}^1(v_1)= \{ v_1, v_2\}$   only if  
   there were  only two elements 
 of $\pi_{i-2}^{i-1} (F)$. As with Case 3a,  
  the generative property of $\pi_{i-2}^{i-1} (v_1)$ for Agent $2$ and 
 either the adjacency to $B_{i-2}$ through agent $2$ 
  or the definition of $w_0$  results in a contradiction.     
   On the other hand, if  $v_1$ is   not generative for Agent 1, by 
  the fact that it is also not proto-generative 
  for Agent 1 (by proximity to $B_{i-1}$ and the definition of $w_0$) 
 $w_1$  would be  equal to $w$, which we assumed 
  to be not in $B_i$.

\vskip.2cm 
\noi   {\bf Case 3c;  $l=1$, meaning 
 that  
  $w$ and $w'$ are extensions of the same $v\in \Omega_{i-1}^f$ :}
 If $v$ were not proto-generative for some  agent then $w$ and $w'$ would 
 already be adjacent. Now assuming that $v$ is generative for 
 both agents, the
 only way for an  extension $w^*$ of $v$ to 
 belong to $B_i$ would be if  $M_{i-1}^j(w^*)= F^j_{i-1}(v)
 \in \overline {\cal F}_{i-1}^j(f)$ for both 
 $j=1,2$, and either   $F^1_{i-1} 
(v)$ 
 or    $F^2_{i-1} (v)$
 contained a   member 
 of $B_{i-1}$ and no member of $A_{i-1}$. If  
  either   $F^1_{i-1} 
(v)$ 
 or    $F^2_{i-1} (v)$
 contained a   member 
 of $B_{i-1}$ and no member of $A_{i-1}$, then by $w$ and $w'$ both not in 
 $B_i$ it must hold that 
 $M_{i-1}^j(w)\not=  F^j_{i-1}(v)$ for some $j$ 
 and  $M_{i-1}^{j'}(w')\not=  F^{j'}_{i-1}(v)$ for some $j'$.
 If  $M_{i-1}^j(w)=  F^j_{i-1}(v)$ 
 and  $M_{i-1}^{j}(w')=  F^{j}_{i-1}(v)$ for the same  $j$ then 
 $w$ and $w'$ were already adjacent. 
  Otherwise if $j$ can be different from $j'$  
   define 
 $w^*\in \Omega^f_i$ by  $M_{i-1}^j(w^*)=  M_{i-1}^j(w)$ 
 and  $M_{i-1}^{j'}(w^*)=  M_{i-1}^{j'}(w')$ and it follows 
 that $w^*$ will connect $w$ and $w'$ without being 
 in $B_i$. \hfill $\Box$
\vskip.2cm

     \noi     Lemma 6  implies that the removal 
         of $B_i$ does not disconnect $\Omega^f_i$. As we will  see, 
 it does not matter that perhaps $\Omega^f_i\backslash B_i$ may 
 be connected through $A_i$. Due to Lemma 1d and Property 4 defining 
 the set $S$, at every level $i$ the set $B_i$ will vastly outnumber 
 the set $A_i$. More importantly,  the 
   extensions  in  $A_{i+1}, A_{i+2}, \dots$ of a member of $A_i$  
 do not involve the least information extensions. 
  Notice that every member of $B_i$ (except 
 for $i=\inf T$) is by 
 definition connected to some member of $A_i$. 
            \vskip.2cm 

\noi 
{\bf Lemma 7:}
If $i\in T$
 and the shortest adjacency paths   
within 
  $\Omega^f_i\backslash B_i$ between
 $w\in \Omega^f_i\backslash \pi_i^ {n_T(i)}(B_{n_T(i)})$ and 
$\pi_i ^ {n_T(i)}(B_{n_T(i)})
\subseteq \Omega^f_i$ are of length $k\geq 1$, 
then there is an $1\leq l\leq k$ with $p^f_{n_T^{l+1}(i)}(w) \in 
B_{n_T^{l+1}(i)}$.    \vskip.2cm

     \noi                            {\bf Proof:}
 We proceed by induction on $k$.  If $k=1$,  let $c\in   
 \pi_i^ {n_T(i)}(B_{n_T(i)})$ be adjacent to $w\not\in 
 \pi_i^{n_T(i)}(B_{n_T(i)})$ and let 
 $j$ be the agent  
  such that $w$ and $c$ share the same member of 
  $\overline {\cal F}^j_i(f)$. 
 The atom  $p^f_{{n_T}(i)}(c)\in B_{{n_T}(i)}$ can not share  
   the same member of $\overline {\cal F}^j_{{n_T}(i)}(f)$ with a member of  
    $A_{{n_T}(i)}$, since by  Lemma 5 the element $c$ would have shared  
     the same member of $\overline {\cal F}^j_i(f)$ with a member of 
      $B_i$ and therefore  
    $p^f_{n_T(i)}(w)$ would also be  a member of $ B_{n_T(i)}$,
 a contradiction to $k=1$. By the definition 
 of $B_{n_T^2(i)}$ it follows that $p^f_{n_T^2(i)} (w) $ is in $B_{n_T^2(i)}$.
  \vskip.2cm 

\noi 
    Assume the claim is true for $k-1\geq 1$.  Let $v\in \Omega^f_i$ 
     be the next element after $c$ in one of the {\bf shortest} 
     adjacency paths  
      within $\Omega^f_i \backslash B_i$ from some $c\in 
    \pi_i^ {n_T(i)}(B_{n_T(i)})$ to $w$. Let $v$ share with $c$ a member of 
      $\overline {\cal F}^j_i(f)$. By the  argument  above    
      we have that $p^f_{ n^2_T(i)}(v) \in B_{n_T^2(i)}$ and 
        also $p^f_{n_T(i)}(v) \not\in B_{n_T(i)}$ and   
        $p^f_{n_T(i)}(u) \not\in B_{n_T(i)}$ for all $u\in \Omega^f_i$ in 
         this adjacency path from $v$ to $w$, including $u=w$ (due 
 to the minimality of the path).   
           Therefore we have an adjacency path of length  
          $k-1$  within $\Omega^f_{n_T(i)} \backslash B_{n_T(i)}$  between 
          $p^f_{n_T(i)}(v)\in \pi_{n_T(i)} \circ 
          \pi ^{-1}_{n_T^2(i)}(B_{n_T^2(i)})$ and  
           $p^f_{n_T(i)}(w)$.  Whether or not it is one 
           of the shortest adjacency 
            paths of this kind we have our conclusion 
            by the induction hypothesis. 
   \hfill $\Box$ \vskip.2cm
\noi Lemma 6 and Lemma 7 show that no member $v$ 
  of $\Omega_i^f$ can be 
 trapped by the formation of the $B_i$,
 meaning that  is impossible 
  for $v$ to belong to $\pi_i (B_{n_T(i)})$
 without $p^f_{n_T(i)} (v)\in B_{n_T(i)}$ 
generating some members of $B_{n_T^2(i)}$. 
  This is critical to the proof of the next lemma. \vskip.2cm 
      \noi {\bf Lemma 8:} The cell 
$C= \{ \gamma (a) \ |\ i\in T, a\in A_i\} $ 
      is  dense
        in ${\bf Ck}(\{ f\})$. \vskip.2cm 
        
\noi         {\bf Proof:} Density of the cell  
 will follow because by choosing any level $i\in T$ and any 
 $w\in \Omega^f_i$ we will  show that 
there is  some level $\hat i\geq i$ with a member of $A_{\hat i}$ 
 extending $w$.
        By Lemma 7 and the fact 
        that $B_i$ does not disconnect  
        $\Omega^f_i \backslash B_i$ (Lemma 6), we need to show that 
       $\pi_i^{n_T(i)}(B_{n_T(i)}) $
        is not empty for every $i\in T$. We establish  this 
 by induction on $i\in T$. The set $B_{\inf T} = \{ w_0\}$ is not 
 empty.  Assume 
 the claim is true for any  particular $i\in T$ and all those members 
 of $T$ before $i$.   
         By Lemma 7 all elements of 
           $ B_i\cup (\pi_i^{n_T(i)}(B_{n_T(i)}))$ are 
           within an 
           adjacency distance of $m:= |T\cap \{ 1,2, \dots , i\}|$ from 
          $\gamma_i (w_0)$, yet the diameter of $\Omega^f_i$ is at least  
          $2m+1$ by the definition of $S$; this means that 
 both $\pi_i^ {n_T(i)}(B_{n_T(i)})$ and  $\Omega^f_i\backslash
 \pi_i^ {n_T(i)}(B_{n_T(i)})$ are not empty.
Also by Lemma 7, the non-emptiness of   
         $\pi_i ^ {n_T(i)}(B_{n_T(i)})$ and the existence of    
        some $w\in 
\Omega^f_i\backslash
 \pi_i^ {n_T(i)}(B_{n_T(i)})$ that is  of positive  
         but finite adjacency distance 
          from $\pi_i ^ {n_T(i)}(B_{n_T(i)})$ implies the 
          non-emptiness 
of $\pi_{n_T(i)}^ { n^2_T(i)}
(B_{ n^2_T(i)})\subseteq 
           \Omega^f _{n_T(i)}$. Density now follows by Lemma 7 and the
 connectivity of all the $\Omega^f_i$.    \hfill $\Box$ 
         \vskip.2cm

    \noi {\bf Theorem 2:}       
           If  the formula $f$ is generative then   
           there is a continuum of cells with finite fanout 
that are dense in ${\bf Ck}(\{ f\})$.
         \vskip.2cm   
         
\noi {\bf Proof:}   
  Let $C$ be a cell created from a subset $T\subseteq S$ with 
 $\inf T = \inf S$. Every possibility set of $C$ is created from 
 a least information extension from some level $i\in T$ applied repetitively 
  to a 
member of $ \pi_i^ {n_T(i)}(B_{n_T(i)})$ (sharing its possibility set 
 in $\Omega^f_i$ with a member of  $B_i$)  
  up to the stage $n_T(i)$,
 followed by a limiting of its size according  to  membership
 in  some possibility set of $\Omega^f_{n_T(i)}$. Due to Property 4 defining 
 $S$ and Lemma 1d, we can read off the subset $T$ from the sizes 
 of the possibility sets in $C$. As there are uncountably many infinite 
 subsets $T$ of $S$ with $\inf T=\inf S$, the theorem is proven.  
\hfill $\Box $ 

\section {Infinite generation and infinite fanout}

Let us review the possibilities for cells of finite fanout.     
 All finite cells are 
 defined by the common knowledge of a single formula (Fagin, Halpern, and 
 Vardi 1991).  Combined with 
  results from Simon (1999, 2001)
 if a cell $C$ has finite fanout it can come in only 
 one of four forms: \newline 
1) $C$ is finite, $F(C)$ is finitely generated and maximal in ${\cal T}$,
 \newline
2) $C$ is infinite, $F(C)$ is finitely generated and not maximal in 
 ${\cal T}$, and $C$ is uncentered, \newline 
3) $C$ is infinite, $F(C)$ is infinitely generated and not maximal 
 in ${\cal T}$, and $C$ is  centered.\newline 
4)  $C$ is infinite, $F(C)$ is infinitely generated, and $C$ is
 uncentered.\newline 
Within Case 4, $F(C)$ may or may not be  maximal, as with the Bernoulli 
 shift space examples presented above.  
How can one distinguish Case 3 from the others? A 
 countable cell is centered if and only if it 
 contains at least one isolated point, (Simon 1999), a straightforward 
 application of Baire Category.  
Difficult  to distinguish 
is Case 2 from Case 4, the distinction being that of finite  vs 
 infinite generation.

\vskip.2cm 
\noi How wild can things be if   finite fanout does not hold? 
Until now, this paper has been concerned with   
  the existence and properties of uncountably many cells that 
 share the same set of formulas in common knowledge.
  This  is qualitatively different   
 from that of uncountably many distinct 
Kripke Structure that map injectively  
 to mutually distinct subsets of  the same 
 cell of $\Omega$.  These    
 issues are  different because  every possibility set of a cell is 
 a compact set, a  property not assumed of the image of a semantic 
 model that maps injectively to $\Omega$.
 This distinction comes into sharp contrast when 
 considering the finite fanout property.  A cell of finite fanout has 
 the {\em surjective} property, meaning that all Kripke Structure 
that map to it must map to it surjectively.  
  A non-surjective cell may 
 offer many  possibilities for disconnected Kripke Structure to map to
 some  cell, but  neither this cell nor these Kripke Structure 
 can   have finite fanout.  This is because the 
 image of a possibility set of 
 a multi-partition with evaluation in $\Omega$ 
must be a dense subset
 of a possibility set of $\Omega$ (Lemma 5, Simon 1999).
\vskip.2cm

\noi   We present a  cell that is centered and  compact, 
 meaning also by Proposition 2 of Simon (1999) 
that it has finite adjacency diameter, and yet there is a Kripke Structure 
 with  uncountably many connected components 
 that maps injectively to  
 this cell. Furthermore, the corresponding set of formulas cannot be 
 finitely generated, since the compactness of the cell $C$ implies 
 the maximality of these formulas in ${\cal T}$ and by Theorem 1 of 
 Simon (2001) this would imply that this cell must be finite.
\vskip.2cm 

\noi To explain our claim, we must first describe Example 3 presented in  
Simon (2001). 
This was an example  of a compact cell homeomorphic to a Cantor set 
  with an adjacency radius of 
 2. To construct this example we let
                     $\Omega$ be $ \Omega (X,\{ 1,2 \})$ 
         and defined  a   sequence  of partitions  
          in the following  way:\newline for every $0<i<\infty$ 
  define  $A_i= \{  p_i (w) 
  \ |\ w\in \Omega_{i-1}\} $. Define ${\cal P}_0 =\{ \Omega\}$  and 
   ${\cal P}_i = {\cal P}_{i-1}\vee \{\pi^{-1}_i (A_i)\ ,\  
   \Omega \backslash \pi_i^{-1} (A_i)\} $.
  We labelled the partitions by 
${\cal B}= ({\cal P}_i\ |\ 0\leq i< \infty)$  and  
 we defined a Kripke Structure $${\cal K}({\cal B})= (\Omega ; 
 ({\cal Q}^j\ |\  j\in \{ 1,2\} ), 
  {\cal P}_{\infty}; X; \overline \psi)$$
 where the partition ${\cal P}_{\infty}$ for the third agent is the limit 
  of the partitions ${\cal P}_i$, (meaning that $z$ and $z'$ share  
  the same member of ${\cal P}_{\infty}$ if and only if 
   they share the same member of ${\cal P}_i$ for every $i<\infty$), 
    with $\overline \psi$  and the ${\cal Q}^j$ for $j=1,2$  the same 
     used to define the 
     Kripke Structure $\Omega$.  
 The third agent can distinguish two  points if and only if  
   the no-information extension was applied  on different 
 stages.
We showed (Simon 2001) that the set $\phi^{{\cal K}({\cal B})}(\Omega)$ is  
 a cell of $\Omega ( X, \{ 1,2,3\})$ equivalent as a  semantic 
 model to ${\cal K}({\cal B})$, and furthermore that 
 the map $\phi^{{\cal K}({\cal B})}$ is a homeomorphism between 
 $\Omega= \Omega (X, \{ 1,2\})$ and the cell that is its image. 
\vskip.2cm

           \noi    
  Define $A:=\{ p^S(w)\ |\ S\in 2_{\infty}^{\bf N_0},  
i\in S, w\in 
 \Omega_i\} \subseteq \Omega = \Omega (X, \{ 1,2\})$,
 the set of  all alienated extensions with 
 respect to the tautologies. 
 Define $B:= \phi^{{\cal K} ({\cal B})}(A)$, 
 the image of the set $A$ in  
  $\Omega (X,\{ 1,2,3\})$. 
We will show that  $B$ defines a Kripke Structure 
 with uncountably many connected components. 
To show this, we need some additional  results from Simon (1999). 
\vskip.2cm 
\noi 
     We define a subset $A\subseteq \Omega$ to be  {\em good}  
       if 
       for every $j\in J$ and every 
       $F\in {\cal Q}^j$ satisfying $F\cap A\not= 
        \emptyset$  it follows that  $F\cap A$ is dense in $F$. 
        By Lemma 5 and Lemma 6 of Simon (1999)  $A$ is good if and only if  
     for every $z\in A\ $ $\phi^{{\cal  V} (A)}(z) =z$ (where the Kripke 
 structure ${\cal V} (A)$ is defined above).  
     
\vskip.2cm 
\noi    
First we show that $B$ is a good subset. 
 Let $z=p^S(w)\in \Omega$ for some $i\in S\in 
 2_{\infty}^{\bf N_0}$ and $w\in \Omega_i$. Let $j\in \{ 1, 2\}$, 
 $z\in F\in {\cal Q}^j= {\cal Q}^j (X, \{ 1,2\})$,
 and $F\cap \pi^{-1}_k (v) \not= \emptyset$ 
 for some $v\in \Omega_k$ with $k\in S$ and $k\geq i$. 
 Since $v$ shares the same member   
  of $\overline {\cal F}^j_k$ with $\pi_k (z)$ we have that $p^S(v) \in F$. 
 Otherwise let $z\in P\in {\cal P}_{\infty}$ and let 
  $P\cap \pi_k^{-1} (v)\not= \emptyset$ 
  for some $v\in \Omega_k$ 
  with $k\in S$ and $k\geq i$.
  Likewise $p^S(v)$ is in $ P$, 
  since $\pi_k^{-1}(v)$ 
  shares the same member of ${\cal P}_k$ with   
   $z=p^S(\pi_k (z))$.         By the above mentioned homeomorphism  
 $B$ is a good subset.

     \vskip.2cm 
\noi                             Fix  $w_0\in \Omega_0 $.
Next we assume  that the adjacency distance 
between $p^{ S}(w_0)$ and $p^{T}(w_0)$
  within   the Kripke Structure ${\cal V}^{{\cal K}({\cal B})} (A)$    
  is $l<\infty$   
  for some pair $S,T\in 2_{\infty}^ {\bf N_0}$   with $S$ and $T$ 
  both containing $\{ 0\}$.  
  Let $p^{S}(w_0)=z_0, z_1, \dots , z_l=p^{T}(w_0)$ be a 
   path of members of $A$ such that for every $0\leq k\leq l-1\ $ 
   $z_k$ and $z_{k+1}$ 
   share the same  
    member of ${\cal Q}^1$, ${\cal Q}^2$, or  
     ${\cal P}_{\infty}$, and 
   for every $0\leq k\leq l\ $  $z_k = p^{S_k}(v_k)$ for  
    $S_k\in 2_{\infty}^{\bf N_0}$ for all $k$,  
    $v_k\in \Omega_{n_k}$ and $n_k \in S_k$ 
    (with  $S_0=S$, $S_l=T$, $n_0=n_l=0$, and $v_0=v_l=w_0$.)  
    Let $N=\max_{0\leq k\leq l}(n_k)$. If   
    $z_k$ and $z_{k+1}$ share the same  
    member of   
     ${\cal P}_{\infty}$ then by the definition of 
       ${\cal P}_i$  we have that  
      $S_k\backslash \{ 0,1, \dots , N-1\} =  
     S_{k+1}\backslash \{ 0,1, \dots , N-1\} $. 
       Now assume that 
    $z_k$ and $z_{k+1}$ share the same  
    member of   
      ${\cal Q}^1$, (respectively ${\cal Q}^2$.) 
If $i\geq \max (n_k, n_{k+1})$ it is not possible for $i$   to be 
 in ${\bf N_0}\backslash S_k$  
   without $i$ being in ${\bf N_0}\backslash S_{k+1}$ (and vice versa).  
        If such an $i\in {\bf N_0}\backslash S_k$ were in $ S_{k+1}$ then 
        $\pi_{i+1}\circ p^{S_{k+1}}(v_{k+1})$ would be 
         a no-information extension and therefore  
         $\pi_{i+1}\circ p^{S_k}(v_k)$ could not share the same member 
          of $\overline {\cal F}^1_{i +1}$ with it, given that $i\not\in 
           S_k$. (We use that all members of $\Omega_i$ for all $i\geq 0$ 
 are generative for both agents.) 
     That suffices for   
      $S_k\backslash \{ 0,1, \dots , N-1\} =  
     S_{k+1}\backslash \{ 0,1, \dots , N-1\} $ for all $k$ 
 and therefore  that    
 $S\backslash \{ 0,1, \dots , N-1\}= 
 T \backslash \{ 0,1, \dots , N-1\}$.  With 
 $\sim$ defined on $2^{\bf N_0}_{\infty}$ as before, 
 we see that $S\not\sim T$ implies that $p^S(w)$ and $p^T(w)$ cannot
  have a finite adjacency distance in the Kripke Structure 
  ${\cal V}^{{\cal K}({\cal B})} (A)$. 

\vskip.2cm 
\noi 
The above argument that $S\backslash \{ 0,1, \dots , N-1\}= 
 T \backslash \{ 0,1, \dots , N-1\}$  works only because 
 all the points  concerned are alienated extensions. 
 With respect to the whole space 
 $\Omega$ the Kripke Structure
 ${\cal K}({\cal B})$ is connected and 
 has an adjacency radius 
 of 2 (see Simon 2001),
 meaning that there is a point such that all other points can be 
 reached from this point by adjacency paths of length 2 or less! 
                  \vskip.2cm

\section {References}
\begin{description}
 \item [Bacharach, M., Gerard-Varet, L.A., Mongin, P, and Shin, H., eds. 
   (1997),] {\em Epistemic Logic and the 
  Theory of Games and Decisions}, Dordrecht, 
   Kluwer.
\item [Breiman, L. (1992),] Probability. Classics in Applied Mathematics, 
 SIAM. 
  \item [Cresswell, M.J. and Hughes, G.E. (1968),] 
 {\em An Introduction to Modal Logic,} Routledge.
 \item [Fagin, R. (1994),] ``A Quantitative Analysis of  Modal Logic," Journal 
  of Symbolic Logic 59, pp. 209-252.
 \item [Fagin, R., Halpern, Y.J. 
and Vardi, M.Y. (1991),] ``A Model-Theoretic 
 Analysis of Knowledge", Journal of the A.C.M. 91 (2), pp. 382-428.
     \item [Halpern, J. and Moses, Y. (1992),] ``A Guide to Completeness and 
     Complexity for Modal Logics of Knowledge and Belief," 
     Artificial Intelligence 54, pp. 319-379.
\item [ Lismont, L. and Mongin, P. (1995),] ``Belief Closure: A Semantics of 
 Common Knowledge for Modal Propositional Logic," Mathematical Social 
  Sciences, Vol. 30, pp. 127-153. 
\item[Mertens, J.-F., Sorin, S.,  Zamir, S. (1994),]  
 {\em Repeated Games},  Core Discussion Papers 9420-22, Universite Catholique 
  de Louvain.
\item [Milgrom, P. and Weber, R. (1985),] ``Distribution Strategies for 
 Games with Incomplete Information", {\it Mathematics of Operations 
 Research} {\bf 10}, pp. 619-632. 
\item [Samet, D. (1990),] ``Ignoring Ignorance and Agreeing to Disagree," 
  Journal of Economic Theory, Vol. 52, No. 1, pp. 190-207.
\item [Simon, R. (1999),] 
``The Difference between Common Knowledge of Formulas and Sets", International
 Journal of Game Theory, Vol. 28, No. 3, pp. 367- 384.
\item [Simon, R. (2001),] 
``The Generation of Formulas in Common Knowledge", International 
 Journal of Game Theory, Vol. 30, pp. 1-18. 
\item [Simon, R. (2003),] ``Games of Incomplete Information, 
Ergodic Theory, and 
the Measurability of  Equilibria", Israel Journal of Mathematics, 
Vol. 138, pp. 73-92.
\item [Sion, M. (1958),] ``On General Minimax Theorems", {\em Pacific 
 Journal of Mathematics}, {\bf 8}, pp. 171-76.

\item [Wagen, S. (1985),]  {\em The Banach-Tarski Paradox}, Cambridge 
 University Press. 

\end {description}

\end{document}